\documentclass[reqno]{article}
\usepackage{amsmath}
\usepackage{color}
\usepackage{float}
\usepackage{cancel}
\providecommand{\MSC}[1]{\textbf{\it{MSC:}} #1}
\providecommand{\keywords}[1]{\textbf{\it{Keywords:}} #1}
\providecommand{\Subject}[1]{\textbf{\it{Subject:}} #1}
\begin{document}

\pagestyle{myheadings}
\newlength{\textlarg}
\newcommand{\barre}[1]{
\settowidth{\textlarg}{#1}
#1\hspace{-\textlarg}%
\rule[0.5ex]{\textlarg}{0.7pt}}

\title{\bf Remarks on the Egyptian 2/D table in favor of a global approach
(D prime number)}
\author{\tt Lionel Br\'{e}hamet \\
Retired research scientist \\
\date{}
{\tt brehamet.l@orange.fr}}
\maketitle
\begin{flushleft}
\begin{abstract}

\noindent
For {\sl h=3} or {\sl 4}, Egyptian decompositions into {\sl h} unit fractions,
like  2/{\sl D} = 1{\sl /D{\scriptsize1}} + ... +1{\sl /D{\scriptsize h}} , were given  by using ({\sl h-1}) divisors ({\sl d{\scriptsize i}}) of 
{\sl D{\scriptsize1}}. This ancient {\sl modus operandi}, well recognized today, provides
{\sl D{\scriptsize i}=DD{\scriptsize1}/d{\scriptsize i}} for {\sl i} greater than {\sl 1}.
Decompositions selected (depending on   {\sl d{\scriptsize i}}) have generally been studied by modern researchers
through the intrinsic  features of  {\sl d{\scriptsize i}} itself. An unconventional method is presented here without considering the {\sl d{\scriptsize i}} properties
but just  the differences  {\sl d{\scriptsize {(h-1)}}}-{\sl d{\scriptsize {h}}}.
In contrast to widespread ideas about the last  denominator like `{\sl D{\scriptsize h}} smaller than 1000', it is more appropriate to adopt a global boundary
of the form `{\sl D{\scriptsize h}} smaller or equal to {\sl 10D}', where {\sl 10} comes from the Egyptian decimal system.
Singular case  {\sl 2/53} (with {\sl15} instead of {\sl10}) is explained. The number of preliminary alternatives before the final decisions is found to be so low 
(71) for {\sl h=3} or {\sl 4} that a detailed overview is possible. A simple additive method of trials, independent of any context, can be carried out, namely
 {\sl 2n+1= d{\scriptsize 2}}+ ... + {\sl d{\scriptsize h}}. Clearly the decisions fit with a minimal value of the differences
 {\sl d{\scriptsize {(h-1)}}}-{\sl d{\scriptsize {h}}}, independently of any  {\sl d{\scriptsize i}} values.	\\ \vspace{1em}

\noindent
\Subject{math.HO}\\
\noindent
\MSC{01A16}\\
\noindent
\end{abstract}

\newpage
									\section*{Preamble}
\hspace*{1.5em}The recto of the Rhind Mathematical Papyrus  (RMP) {\bf  \cite{Peet,Chace,Robins}} contains the so-called Egyptian $2/D$ table. The genesis of a project such as build   
this table will never really be apprehended.  This is not a project as impressive as the construction of a pyramid or temple, however it has been well and truly succeeded. It is impossible  to doubt that pyramid works have not been carried out without a hierarchy of teams well organized in various specialties.
A perfectly organized hierarchy that included team leaders and supervisors. \\
\hspace*{1.5em}It is not hard to imagine that a structured similar organization was also used for the $2/D$ table. This table has not been an exercise in style. It is imperative to keep in mind that it can not be the work of a single scribe, but surely results of indefinite periods of trials and improvements done by an elite team of scribes talented for calculating. 
As it is well known  through dialogues of Plato, the idea of a small number of scholars (philosophers) comes frequently. To these people only, was reserved the right to reflect on issues such as calculations or the study of numbers. He  knew very well that this type of elite was present in the community of scribes of ancient Egypt. He was also aware of their very advanced knowledges  in these areas, but without knowing all secrets. There is no reason today to reject the idea of an elite team or even a chief scribe empowered to decide the last.\\
The time for carrying the table was perhaps over more than a generation
				\footnote{The creative flash of an inspired scholar (ancient or modern) is short. What is generally much longer is the development of the idea 
				and achievement of tools (theoretical or practical) necessary for its application. Of course once the tools lapped their use takes little time!}, 
in order to provide a satisfactory completed product. In such a product nothing should have been left to chance and everything has been deliberately chosen. This is not like a school  exercise where one can use a decomposition rather than another to solve a given problem.\\
\hspace*{1.5em}Once found suitable methods for calculations, it becomes possible to take a look at ``the preliminary draft" in its entirety. This look is necessary  in order to preserve an overall coherence. Some difficulties thus may be highlighted and resolved by a minimum of general decisions, the simplest as possible. 
The number of potential solutions appears as considerably lower than  {\sl ab initio} unrealistic calculations published 
in the  modern literature  {\bf  \cite  {Gillings,BruckSalom}}, namely 22295 or around 28000. We find that it is enough to consider  only $71+71$ possibilities,  then results could be examined before making consistent decisions. This is realistic. A team spirit is very suitable to make obvious the need for a classification and successive resolutions of difficulties encountered during the project progress. Directives given by a leader are implied. 
All these ideas have put us on the track to a comprehensive approach. These ones are the filigree of our analysis.\\

									\section{Data from the papyrus}
\hspace*{1.5em}RMP is also well known by the name of his transcriber, the scribe Ahmes. This latter copied the document around 1650 BCE. The source, now lost, could date from  XIIth dynasty, a golden age of the middle kingdom.
RMP recto  shows  a table of $2$ divided by numbers $D$ from $5$ up to $101$ into "unit fractions". Number $3$ may be considered as implicitly included, because its decomposition is used in the verso for some problems or  it appears elsewhere in Papyrus Kahun {\bf \cite{Imhausen}}. This fact has been commented pertinently by Abdulaziz {\bf \cite{Abdulaziz} }.
\\ \hspace*{1.5em}For $D$ prime only (except number $101$), we present below a reordered excerpt from the $2/D$ table by using {\textcolor{red}{our favorite red numbers $m$}},  that just show the multiplicity of a denominator with $D$.
Please note that {\it they are not the red auxiliary numbers used by Ahmes}, {\sl ie} those ``decoded" by  Gardner {\bf  \cite{GardnerMilo}},
but related with these latter by means of the divisors of the first denominator $D_1$. \\
						\begin{table}[htpb] 
\caption{ \small REORDERED  $2/D$ TABLE FOR PRIME NUMBERS $D$}
\begin{center}
\small
						\begin{tabular}{c}

				\begin{tabular}{ccc}

$\!\!\!\!\!$\begin{tabular}{|c|r|}\hline

$\scriptstyle  2/D=1/D_1+1/D_2 \;\, \sf [2-terms] $ \\ \hline \hline
$2/3=1/2+1/6\,_ {\textcolor{red}{2}}$   \\  \hline
$2/5=1/3+1/15\,_{\textcolor{red}{ 3}}$   \\ \hline
$2/7=1/4+1/28\,_ {\textcolor{red}{4}}$   \\  \hline
$2/11=1/6+1/66\,_{\textcolor{red}{ 6}}$   \\ \hline
$2/23=1/12+1/276\,_{\textcolor{red}{ {12}}}$  \\ \hline 
\end{tabular}
								& 
$\!\!\!\!\!$\begin{tabular}{|c|r}\hline

$\scriptstyle  2/D=1/D_1+1/D_2+1/D_3\;\,  \sf [3-terms] $ \\ \hline \hline
$2/13=1/8+1/52\,_{\textcolor{red}{ 4}}+1/104\,_{\textcolor{red}{ 8}}$   \\ \hline
$2/17=1/12+1/51\,_{\textcolor{red}{ 3}}+1/68\,_{\textcolor{red}{ 4}}$   \\  \hline
$2/19=1/12+1/76\,_{\textcolor{red}{ 4}}+1/114\,_{\textcolor{red}{ 6}}$   \\  \hline
$2/31=1/20+1/124\,_{\textcolor{red}{ 4}}+1/155\,_{\textcolor{red}{ 5}}$   \\ \hline
$2/37=1/24+1/111\,_{\textcolor{red}{ 3}}+1/296\,_{\textcolor{red}{ 8}}$   \\ \hline
$2/41=1/24+1/246\,_{\textcolor{red}{ 6}}+1/328\,_{\textcolor{red}{ 8}}$   \\ \hline
$2/47=1/30+1/141\,_{\textcolor{red}{ 3}}+1/470\,_{\textcolor{red}{10}}$   \\ \hline
$2/53=1/30+1/318\,_{\textcolor{red}{ 6}}+1/795\,_{\textcolor{red}{15}}$   \\  \hline
$2/59=1/36+1/236\,_{\textcolor{red}{ 4}}+1/531\,_{\textcolor{red}{9}}$   \\  \hline
$2/67=1/40+1/335\,_{\textcolor{red}{ 5}}+1/536\,_{\textcolor{red}{8}}$   \\  \hline
$2/71=1/40+1/568\,_{\textcolor{red}{ 8}}+1/710\,_{\textcolor{red}{10}}$   \\  \hline
$2/97=1/56+1/679\,_{\textcolor{red}{ 7}}+1/776\,_{\textcolor{red}{8}}$   \\ \hline 
\end{tabular} 
								&
$\!\!\!\!\!$\begin{tabular}{|c|r|}\hline				

$\scriptstyle 2/D=1/D_1+1/D_2+1/D_3+1/D_4\;\, \sf[4-terms] $ \\ \hline \hline
$2/29=1/24+1/58\,_{\textcolor{red}{ 2}}+1/174\,_{\textcolor{red}{ 6}}+1/232\,_{\textcolor{red}{ 8}}$   \\ \hline
$2/43=1/42+1/86\,_{\textcolor{red}{ 2}}+1/129\,_{\textcolor{red}{ 3}}+1/301\,_{\textcolor{red}{ 7}}$   \\  \hline
$2/61=1/40+1/244\,_{\textcolor{red}{ 4}}+1/488\,_{\textcolor{red}{ 8}}+1/610\,_{\textcolor{red}{ 10}}$   \\ \hline
$2/73=1/60+1/219\,_{\textcolor{red}{ 3}}+1/292\,_{\textcolor{red}{ 4}}+1/365\,_{\textcolor{red}{ 5}}$   \\ \hline
$2/79=1/60+1/237\,_{\textcolor{red}{ 3}}+1/316\,_{\textcolor{red}{ 4}}+1/790\,_{\textcolor{red}{ 10}}$   \\  \hline
$2/83=1/60+1/332\,_{\textcolor{red}{ 4}}+1/415\,_{\textcolor{red}{ 5}}+1/498\,_{\textcolor{red}{ 6}}$   \\ \hline
$2/89=1/60+1/356\,_{\textcolor{red}{ 4}}+1/534\,_{\textcolor{red}{ 6}}+1/890\,_{\textcolor{red}{ 10}}$   \\ \hline
\end{tabular}
				\end{tabular}

							\end{tabular}

\end{center}
\label{Papyrus}
							\end{table}
\pagebreak

					\section { Outlines of a global approach}
\label{PART 0}

Actually the whole $2/D$ project can been viewed as a 3-component set (or 3-phases, if you like).

\sf FIRST: \rm discovery of a unique [2-terms] solution, if $D$ is a prime number.\\
\sf SECOND: \rm for a {\sf sub-project [composite numbers] }from $9$ up to $99$, realize that a mini-table, with just four numbers, enables to derive all the composite numbers by using a  \sl multiplicative \rm operation 
			\footnote{Idea already suggested by Gillings {\bf  \cite{Gillings}}}.  \\
Four numbers, \sf  3, 5, 7, 11 \rm are enough.
For instance  $99$ is reached with  \sf{3}\rm $ \times \mathit {33} $ or
\sf{11}\rm $ \times \mathit {9}$. 

This mini-table, a kind of 'Mother-table', looks as follows: \\

\begin{table}[htp] 
\caption{ \sf Basic Mother-Table}
\begin{center}
\begin{tabular}{|l|}
 \hline
$2/3=1/2+1/6\,_ {\textcolor{red}{2}}$   \\ [0.01in] 
$2/5=1/3+1/15\,_{\textcolor{red}{ 3}}$   \\ [0.01in] 
$2/7=1/4+1/28\,_ {\textcolor{red}{4}}$   \\ [0.01in] 
$\cdots\cdots\cdots\cdots\cdots\cdots\cdot$  \\ [0.01in] 
$2/11=1/6+1/66\,_{\textcolor{red}{ 6}}$   \\ [0.01in] \hline
\end{tabular}
\end{center}

\label{MotherTable}
\end{table}

One sees the first four  two-terms  decompositions of $2/D$. $D$ being prime, the table is \sf unique. \rm\\
In `theory', {\sl except if a better decision should be token,} any fraction $2/D$ ($D$ composite) could be decomposed from this table by dividing a given row
by a convenient number. Consider an example:
$2/65\!=\mbox{\sf [ (row 2 )/ (number 13) ]}=\!1/39+1/195\,_{\textcolor{red}{3 }}  $, what is the solution adopted in the papyrus. As a matter of fact, all decompositions
for the {\sf sub-project} were given in two-terms (except for {\sf 2/95} as a logical consequence  of the guidelines adopted by the scribes, that we will justify properly later)
						\footnote{ All the Egyptian decompositions for composite numbers are analyzed in our second paper {\bf  \cite{Brehamet}}}
.\\
\hspace*{1.5em}As the `Mother-table' has no need to higher value than $11$ for the {\sf sub-project}, we can better understand that, from $13$, it  could have been decided to leave decompositions into  2 terms. \\ 
{\sf THIRD:} nothing does more obstacle to start a main part of the whole project, namely decompositions into 3 (or 4 terms if necessary), 
for all prime numbers starting from 13 until 97.\\
The study carried in this paper is devoted to the third phase.

							\subsection{General presentation}
									\label{PART I}
						
\setcounter{equation}{0}
\renewcommand{\theequation }{\Roman{section}.\arabic{equation}}

\hspace*{1.5em}We could have present the problems in the Egyptian manner, as did Abdulaziz {\bf \cite{Abdulaziz}} like for example
$47 \quad \overline{30}\quad \overline{141}\quad \overline{470}\quad \mbox{what means}\quad 2/47=1/30+1/141+1/470$, 
but we preferred a  modern way, more easily understandable to us today. This is unrelated to the spirit in which we thought.
Consider $D$ as given,  $D_1$ is an unknown value to be found. Assume now that $d_2$, $d_3$, $d_4$ are distinct divisors of  $D_1$, with $d_2> d_3> d_4$.
These are also unknowns to find.\\ 
In order to standardize the notations, $D$ is used for {\sf D}enominators and $d$ for {\sf d}ivisors.\\
Look at the following (modern) equations that {\sf decompose the 'unity'} in $3$  or $4$ parts:

\begin{equation}
\mathbf{1}= \frac{D}{2D_1}+ \frac{d_2}{2D_1}+  \frac{d_3}{2D_1}.
\end{equation}
\begin{equation}
\mathbf{1}= \frac{D}{2D_1}+ \frac{d_2}{2D_1}+  \frac{d_3}{2D_1}+  \frac{d_4}{2D_1}.
\end{equation}
\vspace{0.2em}
It can be viewed under another standpoint like additive operations on integers:
\begin{equation}
2D_1={D}+ d_2+ d_3.
\label{eq:additive3}
\end{equation}
\begin{equation}
2D_1={D}+ d_2+ d_3+ d_4.
\label{eq:additive4}
\end{equation}
Since $d_2$, $d_3$, $d_4$ divide  $D_1$ then we are sure to find Egyptian decompositions. 
Indeed, dividing by $DD_1$ we always get  sums of unit fractions:

\begin{equation}
\frac{2}{D}= \frac{1}{D_1}+ \frac{1}{(D_1/d_2)D}+ \frac{1}{(D_1/d_3)D}.
\label{eq:FEgypt3}
\end{equation}
\begin{equation}
\frac{2}{D}= \frac{1}{D_1}+ \frac{1}{(D_1/d_2)D}+ \frac{1}{(D_1/d_3)D}+ \frac{1}{(D_1/d_4)D}.
\label{eq:FEgypt4}
\end{equation}
This method was apparently followed {\bf \cite {GardnerMilo}} in RMP table for prime numbers $D$ from $13$ up to $97$ .\\
As can be seen, except $D_1$, all denominators of each equation appear as a multiple of $D$, namely 
\begin{equation}
D_i=m_i D, \mbox{\hspace{0.5em}where}\hspace{0.5em}m_i=(D_1/d_i).
\label{eq:Relationmd}
\end{equation}

Let us briefly summarize the possibilities as follows

\begin{equation}
\frac{2}{D}= \frac{1}{D_1}+ \frac{1}{D_2}+ \frac{1}{D_3}.
\label{eq:Egypt3}
\end{equation}
\begin{equation}
\frac{2}{D}= \frac{1}{D_1}+ \frac{1}{D_2}+ \frac{1}{D_3}+ \frac{1}{D_4}.
\label{eq:Egypt4}
\end{equation}
The main task consists in the determination of $D_1$ and the convenient choice of $d_i$, 
from the additive equations (\ref{eq:additive3}) or (\ref{eq:additive4}).
The $d_i$'s are the red auxiliary numbers used by the scribe Ahmes.
\begin{equation}
d_i=\frac{D_1}{m_i}.
\label{eq:Ahmesdi}
\end{equation}

								\section{[2-terms] analysis}
									\label{TwoTerms}
\setcounter{equation}{0}
\begin{equation}
\frac{2}{D}= \frac{1}{D_1}+ \frac{1}{D_2}.
\label{eq:Egypt2}
\end{equation}
The only comment (admiring) on the subject is that the scribes actually 
found the right solution (unique)  to the problem, namely

\begin{equation}
D_1=\frac{D+1}{2}\quad \mbox{and}\quad D_2=\frac{D(D+1)}{2}.
\end{equation}

								\section{[3-terms] analysis}
									\label{ThreeTerms}
\setcounter{equation}{0}
\sf Right now consider the [3-terms] cases. Egyptians gave:\\ \rm
\vspace{1em}
\begin{tabular}{ccll}
\begin{tabular}{|c|}
\hline
\tt Ahmes's selections \rm [3-terms]\\ [0.01in] \hline
$2/13=1/8+1/52\,_{\textcolor{red}{ 4}}+1/104\,_{\textcolor{red}{ 8}}$   \\ [0.01in]  \hline
$2/17=1/12+1/51\,_{\textcolor{red}{ 3}}+1/68\,_{\textcolor{red}{ 4}}$   \\ [0.01in]  \hline
$2/19=1/12+1/76\,_{\textcolor{red}{ 4}}+1/114\,_{\textcolor{red}{ 6}}$   \\ [0.01in]  \hline
$2/31=1/20+1/124\,_{\textcolor{red}{ 4}}+1/155\,_{\textcolor{red}{ 5}}$   \\ [0.01in]  \hline
$2/37=1/24+1/111\,_{\textcolor{red}{ 3}}+1/296\,_{\textcolor{red}{ 8}}$   \\ [0.01in]  \hline
$2/41=1/24+1/246\,_{\textcolor{red}{ 6}}+1/328\,_{\textcolor{red}{ 8}}$   \\ [0.01in]  \hline
$2/47=1/30+1/141\,_{\textcolor{red}{ 3}}+1/470\,_{\textcolor{red}{10}}$   \\ [0.01in]  \hline
$2/53=1/30+1/318\,_{\textcolor{red}{ 6}}+1/795\,_{\textcolor{red}{15}}$   \\ [0.01in]  \hline
$2/59=1/36+1/236\,_{\textcolor{red}{ 4}}+1/531\,_{\textcolor{red}{9}}$   \\ [0.01in]  \hline
$2/67=1/40+1/335\,_{\textcolor{red}{ 5}}+1/536\,_{\textcolor{red}{8}}$   \\ [0.01in]  \hline
$2/71=1/40+1/568\,_{\textcolor{red}{ 8}}+1/710\,_{\textcolor{red}{10}}$   \\ [0.01in]  \hline
$2/97=1/56+1/679\,_{\textcolor{red}{ 7}}+1/776\,_{\textcolor{red}{8}}$   \\ [0.01in]  \hline
\end{tabular}
&
\begin{tabular}{c}
$\Leftarrow$
\end{tabular}
&									
\begin{tabular}{|c|}
\hline
\tt Unity decomposition\\ [0.01in] \hline
$16 = 13 + 2 + 1_{}$   \\ [0.01in]  \hline
$24 = 17 + 4 + 3_{}$   \\ [0.01in]  \hline
$24 = 19 + 3 + 2_{}$   \\ [0.01in]  \hline
$40= 31 + 5 + 4_{}$   \\ [0.01in] \hline
$48 = 37 + 8 + 3_{}$   \\ [0.01in] \hline 
$48 = 41 + 4 + 3_{}$   \\ [0.01in] \hline 
$60 = 47 + 10 + 3_{}$   \\ [0.01in] \hline 
$60 = 53 + 5 + 2_{}$   \\ [0.01in] \hline  
$72 = 59 + 9 + 4_{}$   \\ [0.01in] \hline 
$80 = 67 + 8 + 5_{}$   \\ [0.01in] \hline 
$80 = 71 + 5 + 4_{}$   \\ [0.01in] \hline  
$112 = 97 + 8 + 7_{}$   \\ [0.01in] \hline 

\end{tabular}
					&					
\begin{tabular}{l}
.
\end{tabular}
\end{tabular} \\
\vspace{1em}
\sf The task of finding $D_1$ is rather simple, from the moment when one realizes that it is enough to establish a table of odd numbers $(2n+1)_{|n\geq 1}$ as a sum of two numbers $ d_2 +d_3$, with $d_2>d_3$. This is easy to do and independent of any context.  The table contains $n$ doublets \{$d_2, d_3$\}
and $\sup(d_2)=2n$. One can start with the lowest values as follows: $d_3=1, d_2=2,4,6, \cdots; d_3=2, d_2=3,5,7, \cdots$ and so on.\rm\\

From Eq.(\ref{eq:additive3}) the first candidate possible for $D_1$ starts at an initial value $D_1^0=(D+1 )/2$ as in Fibonnaci's studies {\bf \cite{Fibonacci}}. 
We can search for general solutions of the form
\begin{equation}
D_1^n=D_1^0 + n,
\end{equation}
whence
\begin{equation}
2D_1^n-{D}= 2n+1 =d_2+ d_3.
\label{eq:additive3bis}
\end{equation}
Since one of the two $D_1$ divisors \{$d_2,d_3$\} is even, then $D_1$ can not be odd, it must be even. This was rightly stressed by Bruins {\bf \cite{Bruins}}.
From the first table of doublets, a new table (of trials) is built, where this time doublets are selected if 
$d_2,d_3$ divide $[(D+d_2+d_3)/2]$. This provides a $D_1^n$ possible. In this favorable case, first $D_3$ is calculated by $DD_1/d_3$, then $D_2$ 
by $DD_1/d_2$.\\
For $D$ given, the table of trials defined by the equation just below
\begin{equation}
\mathtt{2n+1=d_2+d_3}, \mbox{\hspace{0.5em} where $d_2$ and $d_3$  divide $D_1^n$ },
\label{eq:dividers3terms}
\end{equation}
is bounded  by a $n_{max}$
							\footnote{It can be proved  that no solution can be found beyond $n =(D-3)/2$.}. 
By simplicity in our tables, $D_1^n$ will not be written as $D_1^n (d_2,d_3)$.\\
Even by hand, a realization  of  this table takes few time.
For example decompositions into 3 terms lead to a total of trials with only  $71$ possibilities! From this low value, it is conceivable to present all results according to an appropriate parameter. Once found a $d_3$, a good idea would be select a $d_2$ the closest as possible of $d_3$. This provides a type of classification never glimpsed to our knowledge. 
Thus, a {\sf key parameter} of our paper is defined as follows:
\begin{equation}
\Delta_{d}= d_2 -d_3.  
\end{equation}
{\sf Remarks: Clearly Eq. (\ref{eq:dividers3terms}) is related to Bruins's method of  ``parts'' redistribution  $d_2,d_3$ {\bf  \cite{Bruins}}.
However our method is {\sl `artisanal'} and does not need to know the arithmetic properties of $D_1$. 
Once $D$ given, $D_1^n$ are  found by trials, without calculations. Unlike to Bruins which sought some forms of $D_1$ for finding then  possible $D$ values. The approach is quite different as well as the reasons justifying the Egyptian choices.}\\  \vspace{0.2em}
\hspace*{2.5em}{\sf Although our conceptual formalism is different from that of Abdulaziz {\bf \cite{Abdulaziz}}, we (fortunately) found some similarities, but also elements without counterpart to us. A welcome unison is  the following:\\
Let us consider its fractional parameter $[R]$ that is crucial for all its analyses. In our notations we find }
\begin{equation}
D_1 [R] = (2D_1 - D) =  2n+1 =d_2+ d_3,
\end{equation}
or equivalently expressed
\begin{equation}
[R] =\frac{1}{(D_1/ d_2)}+\frac{1}{(D_1/ d_3)}.
\end{equation}
{\sf When it is said ``{\sl ... keeping the terms of $[R]$ less than $10$ was an essential part of determining how $2$:$n$ is to be decomposed.}'', this should be  understood as
$(D_1/ d_3)\leq 10$ and formulated for us as
the  condition (\ref{eq:ConditionD1_3}) with a Top-flag $\mbox{\boldmath $\top$ }\!\!_ f^{\;\;[3]}=10$. }(See below for our Top-flag definition)\\
{\sf 
However note that the `necessity' of our Top-flag comes directly from the value of $D$, without constituting  a check on $D_1$. That only follows from  Eq. (\ref {eq:TopFlag3}).\\
\hspace*{1.5em}In contrast, parameter $[Q]$, defined in Ref. {\bf \cite{Abdulaziz}} by $[Q]=1-[R]$, does not appear to us and plays no role in our analyses.
In addition, as the impact of closeness ($\Delta_{d}$) does not seem to have been apprehended, it is clear that our argumentation will generally be different.
Even if, for some 'easy' cases, we agree.}

\hspace*{1.5em}In short,  for producing their final table, we assume that the scribes have analyzed all preliminary trial results before doing their choice among various alternatives, considered in their totality, not individually.  \\

Furthermore, {\sf due to decimal numeration used by ancient Egyptians}, one can easily understand that a boundary with a  Top-flag $\mbox{\boldmath $\top$ }\!\!_ f^{\;\;[3]}$ \rm 
for the last denominator was chosen with a priority value equal to $\mathbf {10}$ 
(if possible according to the results given by trials). \\

\hspace*{1.5em}The idea of a Top-flag is far to be a {\it `deus ex machina'}. It naturally arises if we try to solve the problem of decomposition in full generality. See {\sf Appendix A} for more details.  \\
Chief scribe wisely decided to impose a upper bound  to all  the denominators $D_3$,
such that
\begin{equation}
D_3 \leq D \mbox{\boldmath $\top$ }\!\!_ f^{\;\;[3]}  .
\label{eq:TopFlag3}
\end{equation}
This cut-off beyond  $\mbox{\boldmath $\top$ }\!\!_ f^{\;\;[3]}$ is equivalent to a mathematical condition on $D_1$: 
\begin{equation}
D_1 \leq d_3\,\mbox{\boldmath $\top$ }\!\!_ f^{\;\;[3]} .
\label{eq:ConditionD1_3}
\end{equation}
\hspace*{1.5em}Remark that this condition might be exploited \sf from the beginning \rm of the calculations for avoiding to handle too large denominators $D_3$. 
Simply find $d_3$, find $d_2$, then calculate  $D_1$, if condition ( \ref{eq:ConditionD1_3}) \\
is not fulfilled then quit, do not calculate  $D_3$,  $D_2$ and  go to next values for $d_3$, $d_2$, $D_1$ and so on.\\
Actually, if we follow  the method of trials for finding the good choices in the order $d_3 \rightarrow d_2 \rightarrow D_1$,
we are naturally led to be careful of  the closeness of $d_2$, $d_3$, measured by $\Delta_{d}$. This can suggest the idea of a classification according to increasing values of 
$\Delta_{d}$. \\
Since this classification seriously enlightens many solutions chosen by the scribes, it is not impossible to imagine that this `{\sl artisan method}' was actually followed. This is a plausible hypothesis, valueless of evidence obviously. An advantage is also that a similar classification can be applied to the decompositions into 4 terms  with the same success, 
see Sect. \ref{FourTerms}.\\
The symbol $^{Eg}$ will be used for indicating Egyptian selections in our tables.\\
Let us now display a preliminary table of trials, see Table \ref{Tble3terms71}.\\
\begin{table} [htp]
\caption{\sf Table of trials [3-terms] with increasing order  of $\Delta_{d}$, only 71 possibilities ! }
\begin{center}
\scriptsize
\begin{tabular}{|l|c||l||l||c|l|l|} \hline
\multicolumn{7}{|c|}{\sf Table of trials [3-terms] with increasing order  of $\Delta_{d}$ }\\ \hline
$n$ & $2n+1$ & $d_2$ & $d_3$ & $\textcolor{red}{\Delta_{d}}$ & $D_1^n$ & Possible [3-terms] decompositions \\ [0.01in]  \hline \hline
$1$ & $3$ & $2$ & $1$ & $\mathbf{\textcolor{red}{ 1}}$& $8$  & $ \mathbf {2/13=1/8+1/52\,_{\textcolor{red}{ 4}}+1/104\,_{\textcolor{red}{ 8}}}\;\, ^{Eg}$ \\   \hline 
$1$ & $3$ & $2$ & $1$ & $\mathbf{\textcolor{red}{ 1}}$& $10$ & $ \mathbf {2/17\mathit{_a}=1/10+1/85\,_{\textcolor{red}{ 5}}+1/170\,_{\textcolor{red}{ 10}}}$ \\   \hline
$3$ & $7$ & $4$ & $3$ & $\mathbf{\textcolor{red}{ 1}}$& $12$ & $ \mathbf {2/17\mathit{_b}=1/12+1/51\,_{\textcolor{red}{ 3}}+1/68\,_{\textcolor{red}{ 4}}}\;\, ^{Eg}$ \\   \hline 
$2$ & $5$ & $3$ & $2$ & $\mathbf{\textcolor{red}{ 1}}$& $12$ &   $\mathbf {2/19=1/12+1/76\,_{\textcolor{red}{ 4}}+1/114\,_{\textcolor{red}{ 6}}}\;\, ^{Eg}$ \\   \hline  
$1$ & $3$ & $2$ & $1$ & $\mathbf{\textcolor{red}{ 1}}$& $16$ & $ \mathbf {2/29=1/16+1/232\,_{\textcolor{red}{ 8}}+1/464\,_{\textcolor{red}{ 16}} }$ \\   \hline 
$2$ & $5$ & $3$ & $2$ & $\mathbf{\textcolor{red}{ 1}}$& $18$ & $   \mathbf {2/31\mathit{_a}=1/18+1/186\,_{\textcolor{red}{ 6}}+1/279\,_{\textcolor{red}{ 9}}}$ \\    \hline
$4$ & $9$ & $5$ & $4$ & $\mathbf{\textcolor{red}{ 1}}$& $20$ &  $   \mathbf {2/31\mathit{_b}=1/20+1/124\,_{\textcolor{red}{ 4}}+1/155\,_{\textcolor{red}{ 5}}}\;\, ^{Eg}$ \\   \hline 
$1$ & $3$  & $2$ & $1$ & $\mathbf{\textcolor{red}{ 1}}$& $20$  &  $  \mathbf {2/37=1/20+1/370\,_{\textcolor{red}{ 10}}+1/740\,_{\textcolor{red}{ 20}}}$ \\    \hline 
$1$ & $3$ & $2$ & $1$ & $\mathbf{\textcolor{red}{ 1}}$& $22$ & $   \mathbf {2/41\mathit{_a}=1/22+1/451\,_{\textcolor{red}{ 11}}+1/902\,_{\textcolor{red}{ 22}}}$ \\    \hline
$3$ & $7$ & $4$ & $3$ &  $\mathbf{\textcolor{red}{ 1}}$& $24$ & $   \mathbf {2/41\mathit{_b}=1/24+1/246\,_{\textcolor{red}{ 6}}+1/328\,_{\textcolor{red}{ 8}}}\;\, ^{Eg}$ \\    \hline 
$2$ & $5$ & $3$ & $2$ & $\mathbf{\textcolor{red}{ 1}}$& $24$ & $ \mathbf {2/43=1/24+1/344\,_{\textcolor{red}{ 8}}+1/516\,_{\textcolor{red}{ 12}} }$ \\   \hline 
$1$ & $3$ & $2$ & $1$ & $\mathbf{\textcolor{red}{ 1}}$& $28$ & $   \mathbf {2/53=1/28+1/742\,_{\textcolor{red}{ 14}}+1/1484\,_{\textcolor{red}{ 28}}} $ \\   \hline 
$1$ & $3$ & $2$ & $1$ & $\mathbf{\textcolor{red}{ 1}}$& $32$ & $ \mathbf {2/61=1/32+1/976\,_{\textcolor{red}{ 16}}+1/1952\,_{\textcolor{red}{ 32}} }$ \\   \hline 
$2$ & $5$ & $3$ & $2$ & $\mathbf{\textcolor{red}{ 1}}$& $36$ &  $ \mathbf {2/67=1/36+1/804\,_{\textcolor{red}{ 12}}+1/1206\,_{\textcolor{red}{ 18}}} $\\    \hline 
$4$ & $9$ & $5$ & $4$ & $\mathbf{\textcolor{red}{ 1}}$& $40$ & $ \mathbf {2/71\mathit{_a}=1/40+1/568\,_{\textcolor{red}{ 8}}+1/710\,_{\textcolor{red}{ 10}} }\;\, ^{Eg}$ \\   \hline
$6$ & $13$ & $7$ & $6$ & $\mathbf{\textcolor{red}{ 1}}$& $42$ & $ \mathbf {2/71\mathit{_b}=1/42+1/426\,_{\textcolor{red}{ 6}}+1/497\,_{\textcolor{red}{ 7}}}$ \\   \hline 
$1$ & $3$ & $2$ & $1$ & $\mathbf{\textcolor{red}{ 1}}$& $38$ & $ \mathbf {2/73=1/38+1/1387\,_{\textcolor{red}{ 19}}+1/2274\,_{\textcolor{red}{ 38}} }$ \\   \hline 
$2$ & $5$ & $3$ & $2$ & $\mathbf{\textcolor{red}{ 1}}$& $42$ & $ \mathbf {2/79=1/42+1/1106\,_{\textcolor{red}{ 14}}+1/1659\,_{\textcolor{red}{ 21}} }$ \\   \hline 
$1$ & $3$ & $2$ & $1$ & $\mathbf{\textcolor{red}{ 1}}$& $46$ & $ \mathbf {2/89\mathit{_a}=1/46+1/2047\,_{\textcolor{red}{ 23}}+1/4094\,_{\textcolor{red}{ 46}} }$ \\   \hline
$3$ & $7$ & $4$ & $3$ & $\mathbf{\textcolor{red}{ 1}}$& $48$ & $ \mathbf {2/89\mathit{_b}=1/48+1/1068\,_{\textcolor{red}{ 12}}+1/1424\,_{\textcolor{red}{ 16}} }$ \\   \hline 
$1$ & $3$ & $2$ & $1$ & $\mathbf{\textcolor{red}{ 1}}$& $50$ & $  \mathbf {2/97\mathit{_a}=1/50+1/2425\,_{\textcolor{red}{ 25}}+1/4850\,_{\textcolor{red}{ 50}}} $ \\    \hline
$7$ & $15$ & $8$ & $7$ & $\mathbf{\textcolor{red}{ 1}}$& $56$ & $ \mathbf {2/97\mathit{_b}=1/56+1/679\,_{\textcolor{red}{ 7}}+1/776\,_{\textcolor{red}{ 8}}}\;\, ^{Eg}$ \\    \hline \hline 
$3$ & $7$ & $5$ & $2$ & $\mathbf{\textcolor{red}{ 3}}$& $10$ & $ \mathbf {2/13=1/10+1/26\,_{\textcolor{red}{ 2}}+1/65\,_{\textcolor{red}{ 5}}}$ \\   \hline  
$2$ & $5$ & $4$ & $1$ & $\mathbf{\textcolor{red}{ 3}}$& $12$ &  $ \mathbf {2/19=1/12+1/57\,_{\textcolor{red}{ 3}}+1/228\,_{\textcolor{red}{ 12}}}$ \\    \hline 
$2$ & $5$ & $4$ & $1$ & $\mathbf{\textcolor{red}{ 3}}$& $24$ & $ \mathbf {2/43=1/24+1/258\,_{\textcolor{red}{ 6}}+1/1032\,_{\textcolor{red}{ 24}} }$ \\   \hline 
$3$ & $7$ & $5$ & $2$ & $\mathbf{\textcolor{red}{ 3}}$& $30$ & $   \mathbf {2/53=1/30+1/318\,_{\textcolor{red}{ 6}}+1/795\,_{\textcolor{red}{ 15}}}\;\, ^{Eg}$ \\   \hline 
$2$ & $5$ &  $4$ & $1$ & $\mathbf{\textcolor{red}{ 3}}$& $32$ &  $ \mathbf {2/59=1/32+1/472\,_{\textcolor{red}{ 8}}+1/1888\,_{\textcolor{red}{ 32}}} $ \\    \hline 
$2$ & $5$ & $4$ & $1$ & $\mathbf{\textcolor{red}{ 3}}$& $36$ &  $ \mathbf {2/67\mathit{_a}=1/36+1/603\,_{\textcolor{red}{ 9}}+1/2412\,_{\textcolor{red}{ 36}}} $\\    \hline
$6$ & $13$ & $8$ & $5$ & $\mathbf{\textcolor{red}{ 3}}$& $40$ &  $ \mathbf {2/67\mathit{_b}=1/40+1/335\,_{\textcolor{red}{ 5}}+1/536\,_{\textcolor{red}{ 8}}}\;\, ^{Eg}$ \\    \hline
$3$ & $7$ & $5$ & $2$ & $\mathbf{\textcolor{red}{ 3}}$& $40$ & $ \mathbf {2/73=1/40+1/1584\,_{\textcolor{red}{ 8}}+1/1460\,_{\textcolor{red}{ 20}} }$ \\   \hline 
$2$ & $5$ & $4$ & $1$ & $\mathbf{\textcolor{red}{ 3}}$& $44$ & $ \mathbf {2/83=1/44+1/913\,_{\textcolor{red}{ 11}}+1/3652\,_{\textcolor{red}{ 44}} }$ \\   \hline \hline 
$3$ & $7$ & $6$ & $1$ & $\mathbf{\textcolor{red}{ 5}}$& $12$ & $ \mathbf {2/17=1/12+1/34\,_{\textcolor{red}{ 2}}+1/204\,_{\textcolor{red}{ 12}}}$ \\    \hline
$4$ & $9$ & $7$ & $2$ & $\mathbf{\textcolor{red}{ 5}}$& $14$ &  $ \mathbf {2/19=1/14+1/38\,_{\textcolor{red}{ 2}}+1/133\,_{\textcolor{red}{ 7}}}$ \\   \hline 
$3$ & $7$ & $6$ & $1$ & $\mathbf{\textcolor{red}{ 5}}$& $18$ & $ \mathbf {2/29=1/18+1/87\,_{\textcolor{red}{ 3}}+1/522\,_{\textcolor{red}{ 18}} }$ \\   \hline 
$5$ & $11$ & $8$ & $3$ & $\mathbf{\textcolor{red}{ 5}}$& $24$ & $  \mathbf {2/37=1/24+1/111\,_{\textcolor{red}{ 3}}+1/296\,_{\textcolor{red}{ 8}}}\;\, ^{Eg}$ \\   \hline 
$3$ & $7$ & $6$ & $1$ &  $\mathbf{\textcolor{red}{ 5}}$& $24$ & $   \mathbf {2/41=1/24+1/164\,_{\textcolor{red}{ 4}}+1/984\,_{\textcolor{red}{ 24}}}$ \\    \hline 
$4$ & $9$ & $7$ & $2$ & $\mathbf{\textcolor{red}{ 5}}$& $28$ & $   \mathbf {2/47=1/28+1/188\,_{\textcolor{red}{ 4}}+1/658\,_{\textcolor{red}{ 14}}}$ \\    \hline 
$3$ & $7$ & $6$ & $1$ & $\mathbf{\textcolor{red}{ 5}}$& $30$ & $   \mathbf {2/53=1/30+1/265\,_{\textcolor{red}{ 5}}+1/1590\,_{\textcolor{red}{ 30}}} $ \\   \hline 
$6$ & $13$ & $9$ & $4$ & $\mathbf{\textcolor{red}{ 5}}$& $36$ &  $  \mathbf {2/59=1/36+1/236\,_{\textcolor{red}{ 4}}+1/531\,_{\textcolor{red}{ 9}}}\;\, ^{Eg}$ \\  \hline 
$3$ & $7$ & $6$ & $1$ & $\mathbf{\textcolor{red}{ 5}}$& $48$ & $ \mathbf {2/89=1/48+1/712\,_{\textcolor{red}{ 8}}+1/4272\,_{\textcolor{red}{ 48}} }$ \\   \hline \hline 
$4$ & $9$ & $8$ & $1$ & $\mathbf{\textcolor{red}{ 7}}$& $16$ &  $  \mathbf {2/23=1/16+1/46\,_{\textcolor{red}{ 2}}+1/368\,_{\textcolor{red}{ 16}}}$ \\   \hline 
$6$ & $13$ & $10$ & $3$ & $\mathbf{\textcolor{red}{ 7}}$& $30$ & $   \mathbf {2/47=1/30+1/141\,_{\textcolor{red}{ 3}}+1/470\,_{\textcolor{red}{ 10}}}\;\, ^{Eg}$ \\    \hline 
$5$ & $11$ & $9$ & $2$ & $\mathbf{\textcolor{red}{ 7}}$& $36$ & $ \mathbf {2/61=1/36+1/244\,_{\textcolor{red}{ 4}}+1/1098\,_{\textcolor{red}{ 18}} }$ \\   \hline 
$4$ & $9$ & $8$ & $1$ & $\mathbf{\textcolor{red}{ 7}}$& $40$ & $  \mathbf {2/71=1/40+1/355\,_{\textcolor{red}{ 5}}+1/2840\,_{\textcolor{red}{ 40}}} $ \\   \hline 
$7$ & $15$ & $11$ & $4$ & $\mathbf{\textcolor{red}{ 7}}$& $44$ & $ \mathbf {2/73=1/44+1/292\,_{\textcolor{red}{ 4}}+1/803\,_{\textcolor{red}{ 11}} }$ \\   \hline  
$5$ & $11$ & $9$ & $2$ & $\mathbf{\textcolor{red}{ 7}}$& $54$ & $ \mathbf {2/97=1/54+1/582\,_{\textcolor{red}{ 6}}+1/2619\,_{\textcolor{red}{ 27}}} $ \\  \hline \hline 
$5$ & $11$ & $10$ & $1$ & $\mathbf{\textcolor{red}{ 9}}$& $20$ & $ \mathbf {2/29=1/20+1/58\,_{\textcolor{red}{ 2}}+1/580\,_{\textcolor{red}{ 20}} }$ \\   \hline 
$6$ & $13$ & $11$ & $2$ & $\mathbf{\textcolor{red}{ 9}}$& $22$ & $  \mathbf {2/31=1/22+1/62\,_{\textcolor{red}{ 2}}+1/341\,_{\textcolor{red}{ 11}}}$ \\   \hline 
$5$ & $11$ & $10$ & $1$ & $\mathbf{\textcolor{red}{ 9}}$& $50$ & $ \mathbf {2/89=1/50+1/445\,_{\textcolor{red}{ 5}}+1/4450\,_{\textcolor{red}{ 50}} }$ \\   \hline \hline 
$7$ & $15$ & $13$ & $2$ & $\mathbf{\textcolor{red}{ 11}}$& $26$ & $  \mathbf {2/37=1/26+1/74\,_{\textcolor{red}{ 2}}+1/481\,_{\textcolor{red}{ 13}}}$ \\   \hline 
$6$ & $13$ & $12$ & $1$ & $\mathbf{\textcolor{red}{ 11}}$& $36$ &   $ \mathbf {2/59=1/36+1/177\,_{\textcolor{red}{ 3}}+1/2124\,_{\textcolor{red}{ 36}}} $ \\    \hline 
$8$ & $17$ & $14$ & $3$ & $\mathbf{\textcolor{red}{ 11}}$& $42$ & $  \mathbf {2/67=1/42+1/201\,_{\textcolor{red}{ 3}}+1/938\,_{\textcolor{red}{ 14}}}$ \\    \hline 
$6$ & $13$ & $12$ & $1$ & $\mathbf{\textcolor{red}{ 11}}$& $48$ & $ \mathbf {2/83=1/48+1/332\,_{\textcolor{red}{ 4}}+1/3984\,_{\textcolor{red}{ 48}} }$ \\   \hline 
$7$ & $15$ & $13$ & $2$ & $\mathbf{\textcolor{red}{ 11}}$& $52$ & $ \mathbf {2/89=1/52+1/356\,_{\textcolor{red}{ 4}}+1/2314\,_{\textcolor{red}{ 26}} }$ \\   \hline \hline 
$7$ & $15$ & $14$ & $1$ & $\mathbf{\textcolor{red}{ 13}}$& $28$ & $   \mathbf {2/41=1/28+1/82\,_{\textcolor{red}{ 2}}+1/1148\,_{\textcolor{red}{ 28}}}$ \\   \hline  
$8$ & $17$ & $15$ & $2$ & $\mathbf{\textcolor{red}{ 13}}$& $30$ & $ \mathbf {2/43=1/30+1/86\,_{\textcolor{red}{ 2}}+1/645\,_{\textcolor{red}{ 15}} }$ \\   \hline 
$7$ & $15$ & $14$ & $1$ & $\mathbf{\textcolor{red}{ 13}}$& $56$ & $ \mathbf {2/97=1/56+1/388\,_{\textcolor{red}{ 4}}+1/5432\,_{\textcolor{red}{ 56}}} $ \\    \hline \hline 
$8$ & $17$ & $16$ & $1$ & $\mathbf{\textcolor{red}{ 15}}$& $32$ & $   \mathbf {2/47=1/32+1/94\,_{\textcolor{red}{ 2}}+1/1504\,_{\textcolor{red}{ 32}}} $ \\    \hline 
$8$ & $17$ & $16$ & $1$ & $\mathbf{\textcolor{red}{ 15}}$& $48$ & $ \mathbf {2/79=1/48+1/237\,_{\textcolor{red}{ 3}}+1/3792\,_{\textcolor{red}{ 48}} }$ \\   \hline \hline 
$9$ & $19$ & $18$ & $1$ & $\mathbf{\textcolor{red}{ 17}}$& $36$ & $  \mathbf {2/53=1/36+1/106\,_{\textcolor{red}{ 2}}+1/1908\,_{\textcolor{red}{ 36}}}$ \\    \hline
$9$ & $19$ & $18$ & $1$ & $\mathbf{\textcolor{red}{ 17}}$& $54$ & $ \mathbf {2/89=1/54+1/267\,_{\textcolor{red}{ 3}}+1/4306\,_{\textcolor{red}{ 54}} }$ \\   \hline 
$11$ & $23$ & $20$ & $3$ & $\mathbf{\textcolor{red}{ 17}}$& $60$ & $ \mathbf {2/97=1/60+1/291\,_{\textcolor{red}{ 3}}+1/1940\,_{\textcolor{red}{ 20}}}$ \\   \hline \hline 
$10$ & $21$ & $20$ &$1$ & $\mathbf{\textcolor{red}{ 19}}$& $40$ &  $ \mathbf {2/59=1/40+1/118\,_{\textcolor{red}{ 2}}+1/2360\,_{\textcolor{red}{ 40}}}$  \\    \hline 
$11$ & $23$ & $21$ & $2$ & $\mathbf{\textcolor{red}{ 19}}$& $42$ & $ \mathbf {2/61=1/42+1/122\,_{\textcolor{red}{ 2}}+1/1281\,_{\textcolor{red}{ 21}} }$ \\   \hline \hline 
$12$ & $25$ & $23$ & $2$ & $\mathbf{\textcolor{red}{ 21}}$& $46$ & $  \mathbf {2/67=1/46+1/134\,_{\textcolor{red}{ 2}}+1/1541\,_{\textcolor{red}{ 23}}}$  \\    \hline \hline 
$12$ & $25$ & $24$ & $1$ & $\mathbf{\textcolor{red}{ 23}}$& $48$ & $ \mathbf {2/71=1/48+1/142\,_{\textcolor{red}{ 2}}+1/3408\,_{\textcolor{red}{ 48}}}$ \\    \hline
$13$ & $27$ & $25$ & $2$ & $\mathbf{\textcolor{red}{ 23}}$& $50$ & $ \mathbf {2/73=1/50+1/146\,_{\textcolor{red}{ 2}}+1/1825\,_{\textcolor{red}{ 25}} }$ \\   \hline \hline 
$14$ & $29$ & $27$ & $2$ & $\mathbf{\textcolor{red}{ 25}}$& $54$ & $ \mathbf {2/79=1/54+1/158\,_{\textcolor{red}{ 2}}+1/2133\,_{\textcolor{red}{ 27}} }$ \\   \hline \hline 
$14$ & $29$ & $28$ & $1$ & $\mathbf{\textcolor{red}{ 27}}$& $56$ & $ \mathbf {2/83=1/56+1/166\,_{\textcolor{red}{ 2}}+1/4648\,_{\textcolor{red}{ 56}} }$ \\   \hline \hline 
$15$ & $31$ & $30$ & $1$ & $\mathbf{\textcolor{red}{ 29}}$& $60$ & $ \mathbf {2/89=1/60+1/178\,_{\textcolor{red}{ 2}}+1/5340\,_{\textcolor{red}{ 60}} }$ \\   \hline \hline 
$17$ & $35$ & $33$ & $2$ & $\mathbf{\textcolor{red}{ 31}}$& $66$ & $ \mathbf {2/97=1/66+1/194\,_{\textcolor{red}{ 2}}+1/3201\,_{\textcolor{red}{ 33}}} $ \\    \hline
\end{tabular}
\end{center}

\label{Tble3terms71}
\end{table}
\normalsize \clearpage

As it is clear from Table \ref{Tble3terms71}  an obvious preference for the smallest $\Delta_{d}$ seems to be well followed.  \\ 

After cut-off by $\mbox{\boldmath $\top$ }\!\!_ f^{\;\;[3]}= 10$ Table \ref{Tble3terms71} is reduced and allows us to analyze the following options:
\begin{table} [htbp]
\caption{3-terms options}
\begin{center}
\scriptsize
\begin{tabular}{|l|c||l||l||c|l|l|} \hline
\multicolumn{7}{|c|}{\sf Trials [3-terms] ordered with $\Delta_{d}\nearrow$ showing where are the Egyptian options }\\ \hline
$n$ & $2n+1$ & $d_2$ & $d_3$ & $\textcolor{red}{\Delta_{d}}$ & $D_1^n$ & [3-terms] decompositions $\mathbf {\textcolor{red}{ m_3\leq 10}}$\\ [0.01in]  \hline \hline
$1$ & $3$ & $2$ & $1$ & $\mathbf{\textcolor{red}{ 1}}$& $8$  & $ \mathbf {2/13=1/8+1/52\,_{\textcolor{red}{ 4}}+1/104\,_{\textcolor{red}{ 8}}}\;\, ^{Eg}$ \\   \hline \hline
$1$ & $3$ & $2$ & $1$ & $\mathbf{\textcolor{red}{ 1}}$& $10$ & $ \mathbf {2/17\mathit{_a}=1/10+1/85\,_{\textcolor{red}{ 5}}+1/170\,_{\textcolor{red}{ 10}}}$ \\   \hline
$3$ & $7$ & $4$ & $3$ & $\mathbf{\textcolor{red}{ 1}}$& $12$ & $ \mathbf {2/17\mathit{_b}=1/12+1/51\,_{\textcolor{red}{ 3}}+1/68\,_{\textcolor{red}{ 4}}}\;\, ^{Eg{\textcolor{red}{{\star }}}}$ \\   \hline \hline
$2$ & $5$ & $3$ & $2$ & $\mathbf{\textcolor{red}{ 1}}$& $12$ &   $\mathbf {2/19=1/12+1/76\,_{\textcolor{red}{ 4}}+1/114\,_{\textcolor{red}{ 6}}}\;\, ^{Eg}$ \\   \hline \hline 
$2$ & $5$ & $3$ & $2$ & $\mathbf{\textcolor{red}{ 1}}$& $18$ & $   \mathbf {2/31\mathit{_a}=1/18+1/186\,_{\textcolor{red}{ 6}}+1/279\,_{\textcolor{red}{ 9}}}$ \\    \hline
$4$ & $9$ & $5$ & $4$ & $\mathbf{\textcolor{red}{ 1}}$& $20$ &  $   \mathbf {2/31\mathit{_b}=1/20+1/124\,_{\textcolor{red}{ 4}}+1/155\,_{\textcolor{red}{ 5}}}\;\, ^{Eg{\textcolor{red}{{\star }}}}$ \\   \hline \hline
$3$ & $7$ & $4$ & $3$ &  $\mathbf{\textcolor{red}{ 1}}$& $24$ & $   \mathbf {2/41=1/24+1/246\,_{\textcolor{red}{ 6}}+1/328\,_{\textcolor{red}{ 8}}}\;\, ^{Eg}$ \\    \hline \hline
$4$ & $9$ & $5$ & $4$ & $\mathbf{\textcolor{red}{ 1}}$& $40$ & $ \mathbf {2/71\mathit{_a}=1/40+1/568\,_{\textcolor{red}{ 8}}+1/710\,_{\textcolor{red}{ 10}} }\;\, ^{Eg}$ \\   \hline
$6$ & $13$ & $7$ & $6$ & $\mathbf{\textcolor{red}{ 1}}$& $42$ & $ \mathbf {2/71\mathit{_b}=1/42+1/426\,_{\textcolor{red}{ 6}}+1/497\,_{\textcolor{red}{ 7}}}$ \\   \hline \hline
$7$ & $15$ & $8$ & $7$ & $\mathbf{\textcolor{red}{ 1}}$& $56$ & $ \mathbf {2/97=1/56+1/679\,_{\textcolor{red}{ 7}}+1/776\,_{\textcolor{red}{ 8}}}\;\, ^{Eg{\textcolor{red}{{\star }}}}$ \\    \hline \hline \hline
$3$ & $7$ & $5$ & $2$ & $\mathbf{\textcolor{red}{ 3}}$& $10$ & $ \mathbf {2/13=1/10+1/26\,_{\textcolor{red}{ 2}}+1/65\,_{\textcolor{red}{ 5}}}$ \\   \hline \hline 
$6$ & $13$ & $8$ & $5$ & $\mathbf{\textcolor{red}{ 3}}$& $40$ &  $ \mathbf {2/67=1/40+1/335\,_{\textcolor{red}{ 5}}+1/536\,_{\textcolor{red}{ 8}}}\;\, ^{Eg}$ \\    \hline \hline \hline
$4$ & $9$ & $7$ & $2$ & $\mathbf{\textcolor{red}{ 5}}$& $14$ &  $ \mathbf {2/19=1/14+1/38\,_{\textcolor{red}{ 2}}+1/133\,_{\textcolor{red}{ 7}}}$ \\   \hline \hline
$5$ & $11$ & $8$ & $3$ & $\mathbf{\textcolor{red}{ 5}}$& $24$ & $  \mathbf {2/37=1/24+1/111\,_{\textcolor{red}{ 3}}+1/296\,_{\textcolor{red}{ 8}}}\;\, ^{Eg}$ \\    \hline \hline
$6$ & $13$ & $9$ & $4$ & $\mathbf{\textcolor{red}{ 5}}$& $36$ &  $  \mathbf {2/59=1/36+1/236\,_{\textcolor{red}{ 4}}+1/531\,_{\textcolor{red}{ 9}}}\;\, ^{Eg}$ \\    \hline \hline \hline
$6$ & $13$ & $10$ & $3$ & $\mathbf{\textcolor{red}{ 7}}$& $30$ & $   \mathbf {2/47=1/30+1/141\,_{\textcolor{red}{ 3}}+1/470\,_{\textcolor{red}{ 10}}}\;\, ^{Eg}$ \\    \hline 
\end{tabular}\\ 
\end{center}
\label{3TERMSOPT}
\end{table}\\
This table shows rare instances where multipliers $m_2$, $m_3$  are consecutive. It is always an interesting quality that does not require sophisticated mathematical justification.
That will be  denoted  by a asterisk ${\textcolor{red}{^{\star }}}$. Two instances are found also in [4-terms] series
with $m_2$, $m_3$, $m_4$, see Section \ref{FourTerms}. \\
\nopagebreak [4]
Just as an indication, we display below the cases dropped out of a [3-terms] decomposition:\\ 
\begin{table}[htbp]
\caption{\sf Fractions to be broken down into  4-terms}
\begin{center}
\scriptsize
\begin{tabular}{|l|c||l||l||c|l|l|} \hline
\multicolumn{7}{|c|}{\sf Table of trials [3-terms] for fractions to be broken down into  4-terms}\\ \hline
$n$ & $2n+1$ & $d_2$ & $d_3$ & $\textcolor{red}{\Delta_{d}}$ & $D_1^n$ & Possible [3-terms] decompositions \\   \hline \hline
$4$ & $9$ & $8$ & $1$ & $\mathbf{\textcolor{red}{ 7}}$& $16$ &  $ \mathbf {2/23=1/16+1/46\,_{\textcolor{red}{ 2}}+1/368\,_{\textcolor{red}{ 16}}}$ \\   \hline \hline 
$1$ & $3$ & $2$ & $1$ & $\mathbf{\textcolor{red}{ 1}}$& $16$ & $ \mathbf {2/29=1/16+1/232\,_{\textcolor{red}{ 8}}+1/464\,_{\textcolor{red}{ 16}} }$ \\   \hline
$3$ & $7$ & $6$ & $1$ & $\mathbf{\textcolor{red}{ 5}}$& $18$ & $ \mathbf {2/29=1/18+1/87\,_{\textcolor{red}{ 3}}+1/522\,_{\textcolor{red}{ 18}} }$ \\   \hline
$5$ & $11$ & $10$ & $1$ & $\mathbf{\textcolor{red}{ 9}}$& $20$ & $ \mathbf {2/29=1/20+1/58\,_{\textcolor{red}{ 2}}+1/580\,_{\textcolor{red}{ 20}} }$ \\   \hline \hline
$2$ & $5$ & $4$ & $1$ & $\mathbf{\textcolor{red}{ 3}}$& $24$ & $ \mathbf {2/43=1/24+1/258\,_{\textcolor{red}{ 6}}+1/1032\,_{\textcolor{red}{ 24}} }$ \\   \hline
$2$ & $5$ & $3$ & $2$ & $\mathbf{\textcolor{red}{ 1}}$& $24$ & $ \mathbf {2/43=1/24+1/344\,_{\textcolor{red}{ 8}}+1/516\,_{\textcolor{red}{ 12}} }$ \\   \hline
$8$ & $17$ & $15$ & $2$ & $\mathbf{\textcolor{red}{ 13}}$& $30$ & $ \mathbf {2/43=1/30+1/86\,_{\textcolor{red}{ 2}}+1/645\,_{\textcolor{red}{ 15}} }$ \\   \hline \hline
$1$ & $3$ & $2$ & $1$ & $\mathbf{\textcolor{red}{ 1}}$& $32$ & $ \mathbf {2/61=1/32+1/976\,_{\textcolor{red}{ 16}}+1/1952\,_{\textcolor{red}{ 32}} }$ \\   \hline
$5$ & $11$ & $9$ & $2$ & $\mathbf{\textcolor{red}{ 7}}$& $36$ & $ \mathbf {2/61=1/36+1/244\,_{\textcolor{red}{ 4}}+1/1098\,_{\textcolor{red}{ 18}} }$ \\   \hline
$11$ & $23$ & $21$ & $2$ & $\mathbf{\textcolor{red}{ 19}}$& $42$ & $ \mathbf {2/61=1/42+1/122\,_{\textcolor{red}{ 2}}+1/1281\,_{\textcolor{red}{ 21}} }$ \\   \hline \hline
$1$ & $3$ & $2$ & $1$ & $\mathbf{\textcolor{red}{ 1}}$& $38$ & $ \mathbf {2/73=1/38+1/1387\,_{\textcolor{red}{ 19}}+1/2274\,_{\textcolor{red}{ 38}} }$ \\   \hline
$3$ & $7$ & $5$ & $2$ & $\mathbf{\textcolor{red}{ 3}}$& $40$ & $ \mathbf {2/73=1/40+1/1584\,_{\textcolor{red}{ 8}}+1/1460\,_{\textcolor{red}{ 20}} }$ \\   \hline
$7$ & $15$ & $11$ & $4$ & $\mathbf{\textcolor{red}{ 7}}$& $44$ & $ \mathbf {2/73=1/44+1/292\,_{\textcolor{red}{ 4}}+1/803\,_{\textcolor{red}{ 11}} }$ \\   \hline
$13$ & $27$ & $25$ & $2$ & $\mathbf{\textcolor{red}{ 23}}$& $50$ & $ \mathbf {2/73=1/50+1/146\,_{\textcolor{red}{ 2}}+1/1825\,_{\textcolor{red}{ 25}} }$ \\   \hline \hline
$2$ & $5$ & $3$ & $2$ & $\mathbf{\textcolor{red}{ 1}}$& $42$ & $ \mathbf {2/79=1/42+1/1106\,_{\textcolor{red}{ 14}}+1/1659\,_{\textcolor{red}{ 21}} }$ \\   \hline 
$8$ & $17$ & $16$ & $1$ & $\mathbf{\textcolor{red}{ 15}}$& $48$ & $ \mathbf {2/79=1/48+1/237\,_{\textcolor{red}{ 3}}+1/3792\,_{\textcolor{red}{ 48}} }$ \\   \hline 
$14$ & $29$ & $27$ & $2$ & $\mathbf{\textcolor{red}{ 25}}$& $54$ & $ \mathbf {2/79=1/54+1/158\,_{\textcolor{red}{ 2}}+1/2133\,_{\textcolor{red}{ 27}} }$ \\   \hline \hline
$2$ & $5$ & $4$ & $1$ & $\mathbf{\textcolor{red}{ 3}}$& $44$ & $ \mathbf {2/83=1/44+1/913\,_{\textcolor{red}{ 11}}+1/3652\,_{\textcolor{red}{ 44}} }$ \\   \hline
$6$ & $13$ & $12$ & $1$ & $\mathbf{\textcolor{red}{ 11}}$& $48$ & $ \mathbf {2/83=1/48+1/332\,_{\textcolor{red}{ 4}}+1/3984\,_{\textcolor{red}{ 48}} }$ \\   \hline
$14$ & $29$ & $28$ & $1$ & $\mathbf{\textcolor{red}{ 27}}$& $56$ & $ \mathbf {2/83=1/56+1/166\,_{\textcolor{red}{ 2}}+1/4648\,_{\textcolor{red}{ 56}} }$ \\   \hline \hline
$1$ & $3$ & $2$ & $1$ & $\mathbf{\textcolor{red}{ 1}}$& $46$ & $ \mathbf {2/89=1/46+1/2047\,_{\textcolor{red}{ 23}}+1/4094\,_{\textcolor{red}{ 46}} }$ \\   \hline
$3$ & $7$ & $6$ & $1$ & $\mathbf{\textcolor{red}{ 5}}$& $48$ & $ \mathbf {2/89=1/48+1/712\,_{\textcolor{red}{ 8}}+1/4272\,_{\textcolor{red}{ 48}} }$ \\   \hline
$3$ & $7$ & $4$ & $3$ & $\mathbf{\textcolor{red}{ 1}}$& $48$ & $ \mathbf {2/89=1/48+1/1068\,_{\textcolor{red}{ 12}}+1/1424\,_{\textcolor{red}{ 16}} }$ \\   \hline
$5$ & $11$ & $10$ & $1$ & $\mathbf{\textcolor{red}{ 9}}$& $50$ & $ \mathbf {2/89=1/50+1/445\,_{\textcolor{red}{ 5}}+1/4450\,_{\textcolor{red}{ 50}} }$ \\   \hline
$7$ & $15$ & $13$ & $2$ & $\mathbf{\textcolor{red}{ 11}}$& $52$ & $ \mathbf {2/89=1/52+1/356\,_{\textcolor{red}{ 4}}+1/2314\,_{\textcolor{red}{ 26}} }$ \\   \hline
$9$ & $19$ & $18$ & $1$ & $\mathbf{\textcolor{red}{ 17}}$& $54$ & $ \mathbf {2/89=1/54+1/267\,_{\textcolor{red}{ 3}}+1/4306\,_{\textcolor{red}{ 54}} }$ \\   \hline
$15$ & $31$ & $30$ & $1$ & $\mathbf{\textcolor{red}{ 29}}$& $60$ & $ \mathbf {2/89=1/60+1/178\,_{\textcolor{red}{ 2}}+1/5340\,_{\textcolor{red}{ 60}} }$ \\   \hline 
\end{tabular}
\end{center}

\label{Frac3become4}
\end{table}
\normalsize
\hspace*{1.5em}Our definition of $\mbox{\boldmath $\top$ }\!\!_ f $ does not depend on a arbitrary value of $D_3$ fixed to $1000$
as often assumed in the literature. It depends  only on the circumstances imposed by the current project.  
Subdivide now table  \ref{3TERMSOPT} into 3 sets according to the properties of each $D$. A first with a only one  ${\Delta_{d}}$, a second with  two different ${\Delta_{d}} $ and a third with two conflicting identical ${\Delta_{d}}$.
 That yields:\\
\clearpage
\begin{table} [htbp]
\caption{\sf A single $\Delta_{d}$ [3-terms]\rm }
\begin{center}

\begin{tabular}{|l|c||l||l||c|l|l|} \hline
\multicolumn{6}{|c|}{\sf D with a single $\Delta_{d} $  \rm $\quad$(options: no)} & \multicolumn{1}{l|}{\sf Scribes's decision: obvious}\\ \hline
$n$ & $2n+1$ & $d_2$ & $d_3$ & $\textcolor{red}{\Delta_{d}}$ & $D_1^n$ & $\qquad$ [3-terms] decomposition \\ \hline \hline
$3$ & $7$ & $4$ & $3$ &  $\mathbf{\textcolor{red}{ 1}}$& $24$ & $   \mathbf {2/41=1/24+1/246\,_{\textcolor{red}{ 6}}+1/328\,_{\textcolor{red}{ 8}}}\;\, ^{Eg}$ \\   \hline \hline
$7$ & $15$ & $8$ & $7$ & $\mathbf{\textcolor{red}{ 1}}$& $56$ & $ \mathbf {2/97=1/56+1/679\,_{\textcolor{red}{ 7}}+1/776\,_{\textcolor{red}{ 8}}}\;\, ^{Eg{\textcolor{red}{{\star }}}}$ \\  \hline \hline \hline
$6$ & $13$ & $8$ & $5$ & $\mathbf{\textcolor{red}{ 3}}$& $40$ &  $ \mathbf {2/67=1/40+1/335\,_{\textcolor{red}{ 5}}+1/536\,_{\textcolor{red}{ 8}}}\;\, ^{Eg}$ \\   \hline \hline \hline
$5$ & $11$ & $8$ & $3$ & $\mathbf{\textcolor{red}{ 5}}$& $24$ & $  \mathbf {2/37=1/24+1/111\,_{\textcolor{red}{ 3}}+1/296\,_{\textcolor{red}{ 8}}}\;\, ^{Eg}$ \\    \hline \hline
$6$ & $13$ & $9$ & $4$ & $\mathbf{\textcolor{red}{ 5}}$& $36$ &  $  \mathbf {2/59=1/36+1/236\,_{\textcolor{red}{ 4}}+1/531\,_{\textcolor{red}{ 9}}}\;\, ^{Eg}$ \\   \hline \hline \hline
$6$ & $13$ & $10$ & $3$ & $\mathbf{\textcolor{red}{ 7}}$& $30$ & $   \mathbf {2/47=1/30+1/141\,_{\textcolor{red}{ 3}}+1/470\,_{\textcolor{red}{ 10}}}\;\, ^{Eg}$ \\  \hline 
\end{tabular}\\ 
\end{center}
\label{1Delta3}
\end{table}
\normalsize

\begin{table} [h]
\caption{\sf Two different $\Delta_{d}$ [3-terms]\rm}
\begin{center}
\begin{tabular}{|l|c||l||l||c|l|l|} \hline
\multicolumn{7}{|c|}{\sf D with two different $\Delta_{d}$ \rm $\quad$(options: yes) }\\ \hline
\multicolumn{6}{|c|}{} & \multicolumn{1}{l|}{\sf Scribes's decision: smallest $\Delta_{d}$}\\ \hline
$n$ & $2n+1$ & $d_2$ & $d_3$ & $\textcolor{red}{\Delta_{d}}$ & $D_1^n$ & $\qquad$ [3-terms] decompositions \\   \hline \hline
$1$ & $3$ & $2$ & $1$ & $\mathbf{\textcolor{red}{ 1}}$& $8$  & $ \mathbf {2/13=1/8+1/52\,_{\textcolor{red}{ 4}}+1/104\,_{\textcolor{red}{ 8}}}\;\, ^{Eg}$ \\   \hline \hline
$3$ & $7$ & $5$ & $2$ & $\mathbf{\textcolor{red}{ 3}}$& $10$ & $ \mathbf {2/13=1/10+1/26\,_{\textcolor{red}{ 2}}+1/65\,_{\textcolor{red}{ 5}}}$ \\   \hline \hline 
$2$ & $5$ & $3$ & $2$ & $\mathbf{\textcolor{red}{ 1}}$& $12$ &   $\mathbf {2/19=1/12+1/76\,_{\textcolor{red}{ 4}}+1/114\,_{\textcolor{red}{ 6}}}\;\, ^{Eg}$ \\   \hline \hline 
$4$ & $9$ & $7$ & $2$ & $\mathbf{\textcolor{red}{ 5}}$& $14$ &  $ \mathbf {2/19=1/14+1/38\,_{\textcolor{red}{ 2}}+1/133\,_{\textcolor{red}{ 7}}}$ \\   \hline 
\end{tabular}\\ 
\end{center}
\label{2Delta3}
\end{table}

\begin{table} [h]
\caption{\sf Two conflicting  identical $\Delta_{d}$  [3-terms]\rm}
\begin{center}
\begin{tabular}{|l|c||l||l||c|l|l|} \hline
\multicolumn{7}{|c|}{\sf D with two conflicting identical $\Delta_{d}$ \rm $\quad$(options: yes)}\\ \hline
\multicolumn{6}{|c|}{} & \multicolumn{1}{l|}{\sf Scribes's decision: consecutive multipliers}\\ \hline
$n$ & $2n+1$ & $d_2$ & $d_3$ & $\textcolor{red}{\Delta_{d}}$ & $D_1^n$ & $\qquad$ [3-terms] decompositions \\   \hline \hline
$1$ & $3$ & $2$ & $1$ & $\mathbf{\textcolor{red}{ 1}}$& $10$ & $ \mathbf {2/17\mathit{_a}=1/10+1/85\,_{\textcolor{red}{ 5}}+1/170\,_{\textcolor{red}{ 10}}}$ \\   \hline
$3$ & $7$ & $4$ & $3$ & $\mathbf{\textcolor{red}{ 1}}$& $12$ & $ \mathbf {2/17\mathit{_b}=1/12+1/51\,_{\textcolor{red}{ 3}}+1/68\,_{\textcolor{red}{ 4}}}\;\, ^{Eg{\textcolor{red}{{\star }}}}$ \\   \hline \hline
$2$ & $5$ & $3$ & $2$ & $\mathbf{\textcolor{red}{ 1}}$& $18$ & $   \mathbf {2/31\mathit{_a}=1/18+1/186\,_{\textcolor{red}{ 6}}+1/279\,_{\textcolor{red}{ 9}}}$ \\    \hline
$4$ & $9$ & $5$ & $4$ & $\mathbf{\textcolor{red}{ 1}}$& $20$ &  $   \mathbf {2/31\mathit{_b}=1/20+1/124\,_{\textcolor{red}{ 4}}+1/155\,_{\textcolor{red}{ 5}}}\;\, ^{Eg{\textcolor{red}{{\star }}}}$ \\   \hline 
\end{tabular}\\ \vspace{0.5em}

\begin{tabular}{|l|c||l||l||c|l|l|} \hline
\multicolumn{6}{|c|}{} & \multicolumn{1}{l|}{\sf Scribes's decision: 2n $\leq$ 10}\\ \hline
$n$ & $2n+1$ & $d_2$ & $d_3$ & $\textcolor{red}{\Delta_{d}}$ & $D_1^n$ & $\qquad$ [3-terms] decompositions \\   \hline \hline
$4$ & $9$ & $5$ & $4$ & $\mathbf{\textcolor{red}{ 1}}$& $40$ & $ \mathbf {2/71\mathit{_a}=1/40+1/568\,_{\textcolor{red}{ 8}}+1/710\,_{\textcolor{red}{ 10}} }\;\, ^{Eg}$ \\   \hline
$6$ & $13$ & $7$ & $6$ & $\mathbf{\textcolor{red}{ 1}}$& $42$ & $ \mathbf {2/71\mathit{_b}=1/42+1/426\,_{\textcolor{red}{ 6}}+1/497\,_{\textcolor{red}{ 7}}}
{\;\,^{\textcolor{red}{{\star }}}}$ \\   \hline 
\end{tabular}
\end{center}

\label{11Delta3}
\end{table}
Remark: {\it in the cases involving options possible,  and in these cases only,} \rm the solutions for \\  \{\sf 2/D = 2/13, 2/19, 2/17, 2/31\}  \rm were chosen respectively 
in the set \sf \{$n= 1, 2,  3,  4$\}$_{\mathbf{|}_\mathbf{{\textcolor{red}{2n\leq 10}}}}$.\\\rm
For ruling on  {\sf 2/71} \rm there is  no convincing arithmetical argumentation, then the choice could \\ have been the simplicity and direct observation:   
once again a boundary  like $2n\leq {\textcolor{red}{ 10}}$ is used for picking  $n=4$. That's it. Too simple, but why not? \\ \vspace{0.5em}
\hspace*{1.5em}{\sf After this natural selection by cut-off with a Top-flag $\sf \mbox{\boldmath $\top$ }\!\!_ f^{\;\;[3]}= \sf10$ and appropriate decisions,} 
it remains some cases to be  examined, especially these with $\boxed{10 < m_3 \leq 16}$ because of the singular status of {\sf 2/23,}  
that the scribes will retain with a decomposition into 2 terms. We display below these cases. Of course {\sf 2/61, 2/83} are {\sl ex officio} excluded from the analysis.\\
(Anticipation is made on [4-terms] analysis and related decisions that follow, like $\sf \mbox{\boldmath $\top$ }\!\!_ f^{\;\;[4]}= \sf10$)
\newpage
\begin{table}[htbp]
\caption{Dynamic comparison for transitions $\mathbf {3 \Rightarrow 4}$}
\begin{center}
\scriptsize
\begin{tabular}{|c|r|}
\hline
{\sl Unique} [2-terms] solution \\ [0.01in] \hline
$\mathbf {2/23=1/12+1/276\,_{\textcolor{red}{ {12}}}}\;\, ^{Eg}$   \\ [0.01in] \hline
\end{tabular}\\
\begin{tabular}{|l|c||l||l||c|l|l|} \hline
\multicolumn{7}{|c|}{\sf Selected  trials [3-terms] for $\boxed{ \bf  2/23}\quad {\footnotesize \sl enigma\, ?}\quad \mathbf {(m_3=16)}$}\\ \hline
$n$ & $2n+1$ & $d_2$ & $d_3$ & $\textcolor{red}{\Delta_{d}}$ & $D_1^n$ & {\it Unique} [3-terms] decomposition \\ [0.01in]  \hline \hline
$4$ & $9$ & $8$ & $1$ & $\mathbf{\textcolor{red}{ 7}}$& $16$ &  $ \mathbf {2/23=1/16+\boxed{1/46\,_{\textcolor{red}{ 2}}}+1/368\,_{\textcolor{red}{ 16}}}$ \\   \hline 
\end{tabular}
\scriptsize
\begin{tabular}{|l|c|l||l||l||c|l|l|} \hline
\multicolumn{8}{|c|}{\sf Selected  trials [4-terms] $\boxed{2/23}$}\\ \hline
$n$ & $2n+1$ & $d_2$ & $d_3$ & $d_4$ &$\textcolor{red}{\Delta_{d}^{'}}$ &$D_1^n$ & [4-terms] decomposition $\mathbf {\textcolor{red}{ m_4\leq 10}}$\\ [0.01in]  \hline \hline
$8$ & $17$ & $10$ & $5$ & $2$ & $\mathbf{\textcolor{red}{ 3}}$& $20$ & $\mathbf {2/23=1/20+\boxed{1/46\,_{\textcolor{red}{ 2}}}+1/92\,_{\textcolor{red}{ 4}}+1/230\,_{\textcolor{red}{10}}}$ \\ \hline 
\end{tabular}
\end{center}
\begin{center}
\scriptsize
\begin{tabular}{|l|c||l||l||c|l|l|} \hline
\multicolumn{7}{|c|}{\sf Selected  trials [3-terms] $\boxed{2/29}\qquad \quad \mathbf {(m_3=16)}$}\\ \hline
$n$ & $2n+1$ & $d_2$ & $d_3$ & $\textcolor{red}{\Delta_{d}}$ & $D_1^n$ & Possible [3-terms] decomposition \\ [0.01in]  \hline \hline
$1$ & $3$ & $2$ & $1$ & $\mathbf{\textcolor{red}{ 1}}$& $16$ & $ \mathbf {2/29=1/16+\boxed{1/232\,_{\textcolor{red}{ 8}}}+1/464\,_{\textcolor{red}{ 16}} }$ \\ [0.01in]  \hline
\end{tabular}
\begin{tabular}{|l|c|l||l||l||c|l|l|} \hline
\multicolumn{8}{|c|}{\sf Selected  trials [4-terms] $\boxed{2/29}$}\\ \hline
$n$ & $2n+1$ & $d_2$ & $d_3$ & $d_4$ &$\textcolor{red}{\Delta_{d}^{'}}$ &$D_1^n$ & Possible [4-terms] decompositions $\mathbf {\textcolor{red}{ m_4\leq 10}}$\\ [0.01in]  \hline \hline
$9$ & $19$ & $12$ & $4  $ & $3 $ & $\mathbf{\textcolor{red}{ 1}}$& $24$ & $\mathbf {2/29=1/24+1/58\,_{\textcolor{red}{ 2}}+1/174\,_{\textcolor{red}{ 6}}+
\boxed{1/232\,_{\textcolor{red}{8}}}}\;\, ^{Eg}$ \\ \hline 
$5$ & $11$ & $5$ & $4$ & $2$ & $\mathbf{\textcolor{red}{ 2}}$& $20$ & $ \mathbf {2/29=1/20+1/116\,_{\textcolor{red}{ 4}}+1/145\,_{\textcolor{red}{ 5}}+1/290\,_{\textcolor{red}{ 10}}}$ \\   \hline
$15$ & $31$ & $15$ & $10$ & $6$ & $\mathbf{\textcolor{red}{ 4}}$& $30$ & $ \mathbf {2/29=1/30+1/58\,_{\textcolor{red}{ 2}}+1/87\,_{\textcolor{red}{ 3}}+1/145\,_{\textcolor{red}{5}}}$ \\ \hline 
\end{tabular}
\end{center}
\begin{center}
\scriptsize
\begin{tabular}{|l|c||l||l||c|l|l|} \hline
\multicolumn{7}{|c|}{\sf Selected  trials [3-terms] $\boxed{2/89}\qquad \quad \mathbf {(m_3=16)}$}\\ \hline
$n$ & $2n+1$ & $d_2$ & $d_3$ & $\textcolor{red}{\Delta_{d}}$ & $D_1^n$ & Possible [3-terms] decomposition \\ [0.01in]  \hline \hline
$3$ & $7$ & $4$ & $3$ & $\mathbf{\textcolor{red}{ 1}}$& $48$ & $ \mathbf {2/89\mathit{_b}=1/48+1/1068\,_{\textcolor{red}{ 12}}+1/1424\,_{\textcolor{red}{ 16}} }$ \\   \hline 
\end{tabular}
\begin{tabular}{|l|c|l||l||l||c|l|l|} \hline
\multicolumn{8}{|c|}{\sf Selected  trials [4-terms] $\boxed{2/89}$}\\ \hline
$n$ & $2n+1$ & $d_2$ & $d_3$ & $d_4$ &$\textcolor{red}{\Delta_{d}^{'}}$ &$D_1^n$ & Possible [4-terms] decompositions $\mathbf {\textcolor{red}{ m_4\leq 10}}$\\ [0.01in]  \hline \hline
$15$ & $31$ & $15$ & $10$ & $6$ &  $\mathbf{\textcolor{red}{ 4}}$& $60$ & $ \mathbf {2/89=1/60+1/356\,_{\textcolor{red}{ 4}}+1/534\,_{\textcolor{red}{ 6}}+1/890\,_{\textcolor{red}{10}}}\;\, ^{Eg}$ \\   \hline 
\end{tabular}
\end{center}
\begin{center}
\scriptsize
\begin{tabular}{|l|c||l||l||c|l|l|} \hline
\multicolumn{7}{|c|}{\sf Selected  trials [3-terms] for $\boxed{\bf 2/53}\quad {\footnotesize \sl enigma\, ?}\quad \mathbf {(m_3=15)}$  }\\ \hline
$n$ & $2n+1$ & $d_2$ & $d_3$ & $\textcolor{red}{\Delta_{d}}$ & $D_1^n$ & Possible [3-terms] decomposition \\ [0.01in]  \hline \hline
$3$ & $7$ & $5$ & $2$ & $\mathbf{\textcolor{red}{ 3}}$& $30$ & $   \mathbf {2/53=1/30+\boxed{1/318\,_{\textcolor{red}{ 6}}}+1/795\,_{\textcolor{red}{ 15}}}\;\, ^{Eg}$ \\   \hline
\end{tabular}
\begin{tabular}{|l|c|l||l||l||c|l|l|} \hline
\multicolumn{8}{|c|}{\sf Selected trials [4-terms] for $\boxed{2/53}$ }\\ \hline
$n$ & $2n+1$ & $d_2$ & $d_3$ & $d_4$ &$\textcolor{red}{\Delta_{d}^{'}}$ &$D_1^n$ & [4-terms] decomposition $\mathbf {\textcolor{red}{ m_4\leq 10}}$\\ [0.01in]  \hline \hline
$9$ & $19$ & $9$ & $6  $ & $4 $ & $\mathbf{\textcolor{red}{ 2}}$& $36$ & $\mathbf {2/53=1/36+1/212\,_{\textcolor{red}{ 4}}
+\boxed{1/318\,_{\textcolor{red}{ 6}}}+1/477\,_{\textcolor{red}{9}}}$ \\   \hline 
\end{tabular}
\end{center}
\begin{center}
\scriptsize
\begin{tabular}{|l|c||l||l||c|l|l|} \hline
\multicolumn{7}{|c|}{\sf Selected  trials [3-terms] $\boxed{2/43}\qquad \quad \mathbf {(m_3=12) }$ or $\; \mathbf {(m_3=15) }$}\\ \hline
$n$ & $2n+1$ & $d_2$ & $d_3$ & $\textcolor{red}{\Delta_{d}}$ & $D_1^n$ & Possible [3-terms] decompositions \\ [0.01in]  \hline \hline
$2$ & $5$ & $3$ & $2$ & $\mathbf{\textcolor{red}{ 1}}$& $24$ & $ \mathbf {2/43=1/24+1/344\,_{\textcolor{red}{ 8}}+1/516\,_{\textcolor{red}{ 12}} }$ \\ [0.01in]  \hline
$8$ & $17$ & $15$ & $2$ & \barre{$\mathbf{\textcolor{red}{ 13}}$}& $30$ & \barre{$ \mathbf {2/43=1/30+\boxed{1/86\,_{\textcolor{red}{ 2}}}+1/645\,_{\textcolor{red}{ 15}} }$ }\\ [0.01in]  \hline 
\end{tabular}=
\begin{tabular}{|l|c|l||l||l||c|l|l|} \hline
\multicolumn{8}{|c|}{\sf Selected  trials [4-terms] $\boxed{2/43}$}\\ \hline
$n$ & $2n+1$ & $d_2$ & $d_3$ & $d_4$ &$\textcolor{red}{\Delta_{d}^{'}}$ &$D_1^n$ & [4-terms] decomposition $\mathbf {\textcolor{red}{ m_4\leq 10}}$\\ [0.01in]  \hline \hline
$20$ & $41$ & $21$ & $14$ & $6$ & $\mathbf{\textcolor{red}{ 8}}$& $42$ & $ \mathbf {2/43=1/42+\boxed{1/86\,_{\textcolor{red}{ 2}}}+1/129\,_{\textcolor{red}{ 3}}+1/301\,_{\textcolor{red}{7}}}\;\, ^{Eg}$ \\ \hline 
\end{tabular}
\end{center}

\begin{center}
\scriptsize
\begin{tabular}{|l|c||l||l||c|l|l|} \hline
\multicolumn{7}{|c|}{\sf Selected  trials [3-terms] $\boxed{2/73}\qquad \quad \mathbf {(m_3=11)}$}\\ \hline
$n$ & $2n+1$ & $d_2$ & $d_3$ & $\textcolor{red}{\Delta_{d}}$ & $D_1^n$ & Possible [3-terms] decomposition \\ [0.01in]  \hline \hline
$7$ & $15$ & $11$ & $4$ & $\mathbf{\textcolor{red}{ 7}}$& $44$ & $ \mathbf {2/73=1/44+\boxed{1/292\,_{\textcolor{red}{ 4}}}+1/803\,_{\textcolor{red}{ 11}} }$ \\ [0.01in]  \hline
\end{tabular}
\begin{tabular}{|l|c|l||l||l||c|l|l|} \hline
\multicolumn{8}{|c|}{\sf Selected  trials [4-terms] $\boxed{2/73}$}\\ \hline
$n$ & $2n+1$ & $d_2$ & $d_3$ & $d_4$ &$\textcolor{red}{\Delta_{d}^{'}}$  &$D_1^n$ & [4-terms] decomposition $\mathbf {\textcolor{red}{ m_4\leq 10}}$\\ [0.01in]  \hline \hline
$23$ & $47$ & $20$ & $15$ & $12$ &  $\mathbf{\textcolor{red}{ 3}}$& $60$ & $  \mathbf {2/73\mathit{_c}=1/60+1/219\,_{\textcolor{red}{ 3}}+\boxed{1/292\,_{\textcolor{red}{ 4}}}+1/365\,_{\textcolor{red}{5}}}\;\, ^{Eg{\textcolor{red}{{\star }}}}$ \\ \hline 
\end{tabular}
\end{center}
\end{table}
\normalsize 
\clearpage
We repeat that we are always in a logic of a construction site with difficulties arising in different parts of the project.  
Problems are processed case after case and do not interfere with another previous part. If not, all becomes incomprehensible. A overview supervised by a chief scribe
can not be conflicted. 
The 6 cases presented above confront us with {\sl a dynamic alternative: select the transition from 3 to 4 fractions, or reject it}. This exceptional situation 
is new in the table construction project, as well as the solution itself! It can be observed that 5 cases on 6 have in common the fact that a same denominator 
appears in [3-terms] and [4-terms] decompositions. \\
{\sl A priori}, this fact may be seen as not being an improvement to better decompose a [3-terms] fraction into [4-terms].
Unless we find a real improvement worthwhile.\\
{\sf $\boxed{2/89}$ :} sixth case, out of the category `same denominator', is quickly ruled and [4-terms] decomposition is adopted. (Anyway it belonged to this table only because $m_3=16$). \\
{\sf $\boxed{2/43}$ :} once dropped out the option $m_3=15$, due to a too high gap  $\Delta_{d}=13$,  the same argument holds, then [4-terms] decomposition is adopted.  \\
{\sf $\boxed{2/73}$ :} the [4-terms] expansion provides an improvement since  that leads to three consecutive multipliers \{$3,\,4,\,5$\}, thus this solution is adopted.\\
{\sf Three cases (slightly reordered) remain to be solved,  they are displayed in the following table.}\\

\begin{center}
\scriptsize
\begin{tabular}{|c|r|}
\hline
{\sl Unique} [2-terms] solution \\ [0.01in] \hline
$\mathbf {2/23=1/12+1/276\,_{\textcolor{red}{ {12}}}}\;\, ^{Eg}$   \\ [0.01in] \hline
\end{tabular}\\
\begin{tabular}{|l|c||l||l||c|l|l|} \hline
\multicolumn{7}{|c|}{\sf Selected  trials [3-terms] for $\boxed{ \bf  2/23}\quad {\footnotesize \sl enigma\, ?}\quad \mathbf {(m_3=16)}$}\\ \hline
$n$ & $2n+1$ & $d_2$ & $d_3$ & $\textcolor{red}{\Delta_{d}}$ & $D_1^n$ & {\it Unique} [3-terms] decomposition \\ [0.01in]  \hline \hline
$4$ & $9$ & $8$ & $1$ & $\mathbf{\textcolor{red}{ 7}}$& $16$ &  $ \mathbf {2/23=1/16+\boxed{1/46\,_{\textcolor{red}{ 2}}}+1/368\,_{\textcolor{red}{ 16}}}$ \\   \hline 
\end{tabular}
\scriptsize
\begin{tabular}{|l|c|l||l||l||c|l|l|} \hline
\multicolumn{8}{|c|}{\sf Selected  trials [4-terms] $\boxed{2/23}$}\\ \hline
$n$ & $2n+1$ & $d_2$ & $d_3$ & $d_4$ &$\textcolor{red}{\Delta_{d}^{'}}$ &$D_1^n$ & [4-terms] decomposition $\mathbf {\textcolor{red}{ m_4\leq 10}}$\\ [0.01in]  \hline \hline
$8$ & $17$ & $10$ & $5$ & $2$ & $\mathbf{\textcolor{red}{ 3}}$& $20$ & $\mathbf {2/23=1/20+\boxed{1/46\,_{\textcolor{red}{ 2}}}+1/92\,_{\textcolor{red}{ 4}}+1/230\,_{\textcolor{red}{10}}}$ \\ \hline 
\end{tabular}
\end{center}


\begin{center}
\scriptsize
\begin{tabular}{|l|c||l||l||c|l|l|} \hline
\multicolumn{7}{|c|}{\sf Selected  trials [3-terms] for $\boxed{\bf 2/53}\quad {\footnotesize \sl enigma\, ?}\quad \mathbf {(m_3=15)}$  }\\ \hline
$n$ & $2n+1$ & $d_2$ & $d_3$ & $\textcolor{red}{\Delta_{d}}$ & $D_1^n$ & Possible [3-terms] decomposition \\ [0.01in]  \hline \hline
$3$ & $7$ & $5$ & $2$ & $\mathbf{\textcolor{red}{ 3}}$& $30$ & $   \mathbf {2/53=1/30+\boxed{1/318\,_{\textcolor{red}{ 6}}}+1/795\,_{\textcolor{red}{ 15}}}\;\, ^{Eg}$ \\   \hline
\end{tabular}
\begin{tabular}{|l|c|l||l||l||c|l|l|} \hline
\multicolumn{8}{|c|}{\sf Selected trials [4-terms] for $\boxed{2/53}$ }\\ \hline
$n$ & $2n+1$ & $d_2$ & $d_3$ & $d_4$ &$\textcolor{red}{\Delta_{d}^{'}}$ &$D_1^n$ & [4-terms] decomposition $\mathbf {\textcolor{red}{ m_4\leq 10}}$\\ [0.01in]  \hline \hline
$9$ & $19$ & $9$ & $6  $ & $4 $ & $\mathbf{\textcolor{red}{ 2}}$& $36$ & $\mathbf {2/53=1/36+1/212\,_{\textcolor{red}{ 4}}
+\boxed{1/318\,_{\textcolor{red}{ 6}}}+1/477\,_{\textcolor{red}{9}}}$ \\   \hline 
\end{tabular}
\end{center}

\begin{center}
\scriptsize
\begin{tabular}{|l|c||l||l||c|l|l|} \hline
\multicolumn{7}{|c|}{\sf Selected  trials [3-terms] $\boxed{2/29}\qquad \quad \mathbf {(m_3=16)}$}\\ \hline
$n$ & $2n+1$ & $d_2$ & $d_3$ & $\textcolor{red}{\Delta_{d}}$ & $D_1^n$ & Possible [3-terms] decomposition \\ [0.01in]  \hline \hline
$1$ & $3$ & $2$ & $1$ & $\mathbf{\textcolor{red}{ 1}}$& $16$ & $ \mathbf {2/29=1/16+\boxed{1/232\,_{\textcolor{red}{ 8}}}+1/464\,_{\textcolor{red}{ 16}} }$ \\ [0.01in]  \hline
\end{tabular}
\begin{tabular}{|l|c|l||l||l||c|l|l|} \hline
\multicolumn{8}{|c|}{\sf Selected  trials [4-terms] $\boxed{2/29}$}\\ \hline
$n$ & $2n+1$ & $d_2$ & $d_3$ & $d_4$ &$\textcolor{red}{\Delta_{d}^{'}}$ &$D_1^n$ & Possible [4-terms] decompositions $\mathbf {\textcolor{red}{ m_4\leq 10}}$\\ [0.01in]  \hline \hline
$9$ & $19$ & $12$ & $4  $ & $3 $ & $\mathbf{\textcolor{red}{ 1}}$& $24$ & $\mathbf {2/29\mathit{_a}=1/24+1/58\,_{\textcolor{red}{ 2}}+1/174\,_{\textcolor{red}{ 6}}+
\boxed{1/232\,_{\textcolor{red}{8}}}}\;\, ^{Eg}$ \\ \hline 
$5$ & $11$ & $5$ & $4$ & $2$ & $\mathbf{\textcolor{red}{ 2}}$& $20$ & $ \mathbf {2/29\mathit{_b}=1/20+1/116\,_{\textcolor{red}{ 4}}+1/145\,_{\textcolor{red}{ 5}}+1/290\,_{\textcolor{red}{ 10}}}$ \\   \hline
$15$ & $31$ & $15$ & $10$ & $6$ & $\mathbf{\textcolor{red}{ 4}}$& $30$ & $ \mathbf {2/29\mathit{_c}=1/30+1/58\,_{\textcolor{red}{ 2}}+1/87\,_{\textcolor{red}{ 3}}+1/145\,_{\textcolor{red}{5}}}$ \\ \hline 
\end{tabular}
\end{center}
For each fraction the same denominators (inside a box) have a well defined position in a [3-terms] expansion and another in a [4-terms]. We denote respectively these positions by 
 \sf rank$^{[3]}$ {\rm and} rank$^{[4]}$.\rm \\
Same denominators will be denoted by $\boxed{same D_i}$. The table below summarizes the situation.\\ \vspace{0.5em}
\begin{center}
\begin{tabular} {|c|c|c|c|l|}\hline
\sf Fraction & $\boxed{same D_i}$ & \sf rank$^{[3]}$ &  \sf rank$^{[4]}$ & Appreciation on ranks\\ \hline
\sf 2/23 & $\boxed{46}$ & $\mathbf 2 $ & $\mathbf 2 $ & no interest \\ \hline
\sf 2/53 & $\boxed{318}$ & $\mathbf 2 $ & $\mathbf 3 $ & too near \\ \hline
\sf 2/29$\mathit{_a} $& $\boxed{232}$ & $\mathbf 2 $ & $\mathbf 4 $ & acceptable + smallest $\textcolor{red}{\Delta_{d}^{'}}$ \\ \hline
\end{tabular}	\\
\end{center}
Some convenient rulings ensue, namely \\
{\sf 2/23}; no solution; then come back to the only one solution in 2 terms.\\
{\sf 2/53}; maintain [3-terms] solution; reject [4-terms] solution. \\
{\sf 2/29$\mathit{_a} $}; adopt [4-terms] solution. 
\newpage
 
							\section{[4-terms] analysis}
									\label{FourTerms}
\setcounter{equation}{0}
Right now consider the [4-terms] cases. Egyptians gave:\\  
\vspace{1.5em}
\begin{tabular}{ccll}
\begin{tabular}{|c|}
\hline
\tt Ahmes's selections \rm [4-terms]\\ [0.01in] \hline
$2/29=1/24+1/58\,_{\textcolor{red}{ 2}}+1/174\,_{\textcolor{red}{ 6}}+1/232\,_{\textcolor{red}{ 8}}$   \\ [0.01in] \hline
$2/43=1/42+1/86\,_{\textcolor{red}{ 2}}+1/129\,_{\textcolor{red}{ 3}}+1/301\,_{\textcolor{red}{ 7}}$   \\ [0.01in] \hline
$2/61=1/40+1/244\,_{\textcolor{red}{ 4}}+1/488\,_{\textcolor{red}{ 8}}+1/610\,_{\textcolor{red}{ 10}}$   \\ [0.01in] \hline
$2/73=1/60+1/219\,_{\textcolor{red}{ 3}}+1/292\,_{\textcolor{red}{ 4}}+1/365\,_{\textcolor{red}{ 5}}$   \\ [0.01in] \hline
$2/79=1/60+1/237\,_{\textcolor{red}{ 3}}+1/316\,_{\textcolor{red}{ 4}}+1/790\,_{\textcolor{red}{ 10}}$   \\ [0.01in] \hline
$2/83=1/60+1/332\,_{\textcolor{red}{ 4}}+1/415\,_{\textcolor{red}{ 5}}+1/498\,_{\textcolor{red}{ 6}}$   \\ [0.01in] \hline
$2/89=1/60+1/356\,_{\textcolor{red}{ 4}}+1/534\,_{\textcolor{red}{ 6}}+1/890\,_{\textcolor{red}{ 10}}$   \\ [0.01in] \hline
\end{tabular}
&
\begin{tabular}{c}
$\Leftarrow$
\end{tabular}
&									
\begin{tabular}{|c|}
\hline
\tt Unity decomposition\\ [0.01in] \hline
$48 = 29 + 12 + 4 + 3_{}$   \\ [0.01in]  \hline
$84 = 43 + 21 + 14 + 6_{}$   \\ [0.01in]  \hline
$80 = 61 + 10 + 5 + 4_{}$   \\ [0.01in]  \hline
$120 = 73 + 20 + 15 +12_{}$   \\ [0.01in] \hline 
$120 = 79 + 20 + 15 +6_{}$   \\ [0.01in] \hline 
$120 = 83 + 15 + 12 + 10_{}$   \\ [0.01in] \hline 
$120 = 89 + 15 + 10 + 6_{}$   \\ [0.01in] \hline 
\end{tabular}
					&					\begin{tabular}{l}
.
\end{tabular}
\end{tabular} \\
\vspace{1.5em}
{\sf The task of finding $D_1$ is rather simple, from the moment when one realizes that it is enough to establish a table of odd numbers $(2n+1)_{|n\geq 3}$ as a sum of three numbers $ d_2 +d_3+d_4$, with $d_2>d_3>d_4$. This is easy to do and independent of any context.  
The table contains ($\boldsymbol[\frac{n}{2}\boldsymbol]\boldsymbol[\frac{n+1}{2}\boldsymbol]$ \!-$1$) triplets \{$d_2, d_3,d_4$\}
and $\sup(d_2)=2n$-$2$.} Square brackets here $ \boldsymbol [\;\boldsymbol]$ means `integral part of'. {\sf One can start with the lowest values as follows:
 $d_4=1, d_3=2,3,4, \cdots, d_2=3,4,5, \cdots; d_4=2, d_3=3,4,5, \cdots, d_2=4,5,6, \cdots$ and so on, with the condition $d_3+d_2 \equiv d_4+1 \mod(2)$. }\\
From Eq.(\ref{eq:additive4}) the first candidate possible for $D_1$ starts at the value $D_1^0=(D+1 )/2$. We can search for general solutions of the form
\begin{equation}
D_1^n=D_1^0 + n,
\end{equation}
whence
\begin{equation}
2D_1^n-{D}= 2n+1 =d_2+ d_3+d_4.
\label{eq:additive4bis}
\end{equation}
From the first table of triplets, a new table (of trials) is built, where this time triplets are selected if 
$d_2,d_3,d_4$ divide $[(D+d_2+d_3+d_4)/2]$. This provides a $D_1^n$ possible. In this favorable case, first $D_4$ is calculated by $DD_1/d_4$,
then $D_3$  by $DD_1/d_3$, and $D_2$  by $DD_1/d_2$.\\
This table of trials, properly defined by the equation just below (included the constraints), ie
\begin{equation}
\mathtt{2n+1=d_2+d_3+d_4}, \mbox{\hspace{0.5em} where $d_2$,  $d_3$ and $d_4$  divide $D_1^n$ },
\label{eq:dividers4terms}
\end{equation}
is obviously a bit longer to establish than for doublets.
By simplicity $D_1^n$ will be not written as $D_1^n (d_2,d_3,d_4)$. 
For decompositions into 4 terms the total of trials yields only  $71$ possibilities ! \\
Of course our remark previously made about doublets is still valid for triplets. Likewise, Abdulaziz's parameter [R] takes the form
\begin{equation}
[R] =\frac{1}{(D_1/ d_2)}+\frac{1}{(D_1/ d_3)}+\frac{1}{(D_1/ d_4)}.
\end{equation}
The notation used in our tables will be
\begin{equation}
\Delta_{d}^{'}= d_3 -d_4,
\end{equation}

Chief scribe wisely decided to impose a upper bound  to all  the denominators $D_4$,
such that \rm
\begin{equation}
D_4 \leq D \mbox{\boldmath $\top$ }\!\!_ f^{\;\;[4]}  .
\end{equation}
This cut-off beyond  $\mbox{\boldmath $\top$ }\!\!_ f^{\;\;[4]}$ is equivalent to a mathematical condition on $D_1$: 
\begin{equation}
D_1 \leq d_4\,\mbox{\boldmath $\top$ }\!\!_ f^{\;\;[4]} .
\label{eq:ConditionD1_4}
\end{equation}
Here again, choosing $\mbox{\boldmath $\top$ }\!\!_ f^{\;\;[4]}=10$  is quite appropriate.
Thus a general coherence is ensured throughout the project, since 11 out of 12 decompositions into 3 terms were solved with 
$\mbox{\boldmath $\top$ }\!\!_ f^{\;\;[3]}=10$.\\
\hspace*{1.5em}Remark that the condition (\ref{eq:ConditionD1_4})  might be exploited \sf from the beginning \rm of the calculations for avoiding to handle too large denominators $D_4$. Simply find
$d_4$, find $d_3$, find $d_2$, calculate $D_1$,  
if (\ref{eq:ConditionD1_4}) is not fulfilled then quit, do not calculate  $D_4$,  $D_3$, $D_2$ and  go to next values for $d_4$, $d_3$, $d_2$, $D_1$ etc.

\begin{table}[htp]
\caption{\sf Table of trials [4-terms] with increasing order  of $\Delta_{d}^{'}$, only 71 possibilities ! }
\scriptsize
\begin{center}
\begin{tabular}{|l|c|l||l||l||c|l|l|} \hline
\multicolumn{8}{|c|}{\sf Trials [4-terms] with increasing order of $\Delta_{d}^{'}$ }\\ \hline
$n$ & $2n+1$ & $d_2$ & $d_3$ & $d_4$ &$\textcolor{red}{\Delta_{d}^{'}}$ &$D_1^n$ & Possible [4-terms] decompositions \\ [0.01in]  \hline \hline
$9$ & $19$ & $12$ & $4  $ & $3 $ & $\mathbf{\textcolor{red}{ 1}}$& $24$ & $\mathbf {2/29=1/24+1/58\,_{\textcolor{red}{ 2}}+1/174\,_{\textcolor{red}{ 6}}+1/232\,_{\textcolor{red}{8}}}\;\, ^{Eg}$ \\ \hline 
$5$ & $11$ & $6$ & $3$ & $2$ &  $\mathbf{\textcolor{red}{ 1}}$& $36$ & $\mathbf {2/61\mathit{_a}=1/36+1/366\,_{\textcolor{red}{ 6}}+1/732\,_{\textcolor{red}{ 12}}+1/1098\,_{\textcolor{red}{ 18}}}$ \\ \hline
$9$ & $19$ & $10$ & $5$ & $4$ &  $\mathbf{\textcolor{red}{ 1}}$& $40$ & $\mathbf {2/61\mathit{_b}=1/40+1/244\,_{\textcolor{red}{ 4}}+1/488\,_{\textcolor{red}{ 8}}+1/610\,_{\textcolor{red}{ 10}}}\;\, ^{Eg}$ \\ \hline 
$3$ & $7$ & $4$ & $2$ & $1$ &  $\mathbf{\textcolor{red}{ 1}}$& $40$ & $ \mathbf {2/73\mathit{_a}=1/40+1/730\,_{\textcolor{red}{ 10}}+1/1460\,_{\textcolor{red}{ 20}}+1/2920\,_{\textcolor{red}{ 40}}}$ \\ \hline
$5$ & $11$ & $6$ & $3$ & $2$ &  $\mathbf{\textcolor{red}{ 1}}$& $42$ & $ \mathbf {2/73\mathit{_b}=1/42+1/511\,_{\textcolor{red}{ 7}}+1/1022\,_{\textcolor{red}{ 14}}+1/1533\,_{\textcolor{red}{ 21}}}$ \\ \hline
$11$ & $23$ & $16$ & $4$ & $3$ &  $\mathbf{\textcolor{red}{ 1}}$& $48$ & $ \mathbf {2/73\mathit{_c}=1/48+1/219\,_{\textcolor{red}{ 3}}+1/876\,_{\textcolor{red}{ 12}}+1/1168\,_{\textcolor{red}{ 16}}}$ \\ \hline
$8$ & $17$ & $12$ & $3$ & $2$ & $\mathbf{\textcolor{red}{ 1}}$& $48$ & $ \mathbf {2/79\mathit{_a}=1/48+1/316\,_{\textcolor{red}{ 4}}+1/1264\,_{\textcolor{red}{ 16}}+1/1896\,_{\textcolor{red}{ 24}}}$ \\ \hline
$20$ & $41$ & $30$ & $6$ & $5$ & $\mathbf{\textcolor{red}{ 1}}$& $60$ & $ \mathbf {2/79\mathit{_b}=1/60+1/158\,_{\textcolor{red}{ 2}}+1/790\,_{\textcolor{red}{ 10}}+1/948\,_{\textcolor{red}{ 12}}}$ \\ \hline 
$6$ & $13$ & $8$ & $3$ & $2$ & $\mathbf{\textcolor{red}{ 1}}$& $48$ & $ \mathbf {2/83\mathit{_a}=1/48+1/498\,_{\textcolor{red}{ 6}}+1/1328\,_{\textcolor{red}{ 16}}+1/1992\,_{\textcolor{red}{ 24}}}$ \\ \hline
$6$ & $13$ & $6$ & $4$ & $3$ & $\mathbf{\textcolor{red}{ 1}}$& $48$ & $ \mathbf {2/83\mathit{_b}=1/48+1/664\,_{\textcolor{red}{ 8}}+1/996\,_{\textcolor{red}{ 12}}+1/1328\,_{\textcolor{red}{ 16}}}$ \\ \hline
$14$ & $29$ & $14$ & $8$ & $7$ & $\mathbf{\textcolor{red}{ 1}}$& $56$ & $ \mathbf {2/83\mathit{_c}=1/56+1/332\,_{\textcolor{red}{ 4}}+1/581\,_{\textcolor{red}{ 7}}+1/664\,_{\textcolor{red}{ 8}}}$ \\ \hline 
$18$ & $37$ & $30$ & $4$ & $3$ & $\mathbf{\textcolor{red}{ 1}}$& $60$ & $ \mathbf {2/83\mathit{_d}=1/60+1/166\,_{\textcolor{red}{ 2}}+1/1245\,_{\textcolor{red}{ 15}}+1/1660\,_{\textcolor{red}{20}}}$ \\ \hline
$3$ & $7$ & $4$ & $2$ & $1$ &  $\mathbf{\textcolor{red}{ 1}}$& $48$ & $ \mathbf {2/89\mathit{_a}=1/48+1/1068\,_{\textcolor{red}{ 12}}+1/2136\,_{\textcolor{red}{ 24}}+1/4272\,_{\textcolor{red}{ 48}}}$ \\ \hline
$15$ & $31$ & $20$ & $6$ & $5$ &  $\mathbf{\textcolor{red}{ 1}}$& $60$ & $ \mathbf {2/89\mathit{_b}=1/60+1/267\,_{\textcolor{red}{ 3}}+1/890\,_{\textcolor{red}{ 5}}+1/1068\,_{\textcolor{red}{12}}}$ \\ \hline \hline 
$6$ & $13$ & $9$ & $3 $ & $1$ & $\mathbf{\textcolor{red}{ 2}}$& $18$ & $\mathbf {2/23=1/18+1/46\,_{\textcolor{red}{ 2}}+1/138\,_{\textcolor{red}{ 6}}+1/414\,_{\textcolor{red}{18}}}$ \\ \hline 
$5$ & $11$ & $5$ & $4$ & $2$ & $\mathbf{\textcolor{red}{ 2}}$& $20$ & $ \mathbf {2/29=1/20+1/116\,_{\textcolor{red}{ 4}}+1/145\,_{\textcolor{red}{ 5}}+1/290\,_{\textcolor{red}{ 10}}}$ \\ \hline 
$6$ & $13$ & $7$ & $4$ & $2$ & $\mathbf{\textcolor{red}{ 2}}$& $28$ & $ \mathbf {2/43=1/28+1/172\,_{\textcolor{red}{ 4}}+1/301\,_{\textcolor{red}{ 7}}+1/602\,_{\textcolor{red}{ 14}}}$ \\ \hline 
$8$ & $17$ & $13$ & $3$ & $1$ &  $\mathbf{\textcolor{red}{ 2}}$& $39$ & $\mathbf {2/61=1/39+1/183\,_{\textcolor{red}{ 3}}+1/793\,_{\textcolor{red}{ 13}}+1/2379\,_{\textcolor{red}{ 39}}}$ \\ \hline
$5$ & $11$ & $7$ & $3$ & $1$ &  $\mathbf{\textcolor{red}{ 2}}$& $42$ & $ \mathbf {2/73\mathit{_a}=1/42+1/438\,_{\textcolor{red}{ 6}}+1/1022\,_{\textcolor{red}{ 14}}+1/3066\,_{\textcolor{red}{ 42}}}$ \\ \hline
$8$ & $17$ & $9$ & $5$ & $3$ &  $\mathbf{\textcolor{red}{ 2}}$& $45$ & $ \mathbf {2/73\mathit{_b}=1/45+1/365\,_{\textcolor{red}{ 5}}+1/657\,_{\textcolor{red}{ 9}}+1/1095\,_{\textcolor{red}{ 15}}}$ \\ \hline 
$18$ & $37$ & $15$ & $12$ & $10$ & $\mathbf{\textcolor{red}{ 2}}$& $60$ & $ \mathbf {2/83=1/60+1/332\,_{\textcolor{red}{ 4}}+1/415\,_{\textcolor{red}{ 5}}+1/498\,_{\textcolor{red}{6}}}\;\, ^{Eg}$ \\ \hline 
$18$ & $37$ & $21$ & $9$ & $7$ &  $\mathbf{\textcolor{red}{ 2}}$& $63$ & $ \mathbf {2/89=1/63+1/267\,_{\textcolor{red}{ 3}}+1/623\,_{\textcolor{red}{ 7}}+1/801\,_{\textcolor{red}{9}}}$ \\ \hline \hline 
$8$ & $17$ & $10$ & $5$ & $2$ & $\mathbf{\textcolor{red}{ 3}}$& $20$ & $\mathbf {2/23=1/20+1/46\,_{\textcolor{red}{ 2}}+1/92\,_{\textcolor{red}{ 4}}+1/230\,_{\textcolor{red}{10}}}$ \\ \hline 
$8$ & $17$ & $10$ & $5$ & $2$ & $\mathbf{\textcolor{red}{ 3}}$& $30$ & $ \mathbf {2/43\mathit{_a}=1/30+1/129\,_{\textcolor{red}{ 3}}+1/258\,_{\textcolor{red}{ 6}}+1/645\,_{\textcolor{red}{15}}}$ \\ \hline
$10$ & $21$ & $16$ & $4$ & $1$ & $\mathbf{\textcolor{red}{ 3}}$& $32$ & $ \mathbf {2/43\mathit{_b}=1/32+1/86\,_{\textcolor{red}{ 2}}+1/344\,_{\textcolor{red}{ 8}}+1/1376\,_{\textcolor{red}{32}}}$ \\ \hline 
$5$ & $11$ & $6$ & $4$ & $1$ &  $\mathbf{\textcolor{red}{ 3}}$& $36$ & $\mathbf {2/61\mathit{_a}=1/36+1/366\,_{\textcolor{red}{ 6}}+1/549\,_{\textcolor{red}{ 9}}+1/2196\,_{\textcolor{red}{ 36}}}$ \\ \hline
$11$ & $23$ & $14$ & $6$ & $3$ &  $\mathbf{\textcolor{red}{ 3}}$& $42$ & $ \mathbf {2/61\mathit{_b}=1/42+1/183\,_{\textcolor{red}{ 3}}+1/427\,_{\textcolor{red}{ 7}}+1/854\,_{\textcolor{red}{14}}}$ \\ \hline 
$13$ & $27$ & $22$ & $4$ & $1$ &  $\mathbf{\textcolor{red}{ 3}}$& $44$ & $ \mathbf {2/61\mathit{_c}=1/44+1/122\,_{\textcolor{red}{ 2}}+1/671\,_{\textcolor{red}{ 11}}+1/2684\,_{\textcolor{red}{44}}}$ \\ \hline
$15$ & $31$ & $26$ & $4$ & $1$ &  $\mathbf{\textcolor{red}{ 3}}$& $52$ & $ \mathbf {2/73\mathit{_a}=1/52+1/146\,_{\textcolor{red}{ 2}}+1/949\,_{\textcolor{red}{ 13}}+1/3796\,_{\textcolor{red}{ 52}}}$ \\ \hline
$19$ & $39$ & $28$ & $7$ & $4$ &  $\mathbf{\textcolor{red}{ 3}}$& $56$ & $ \mathbf {2/73\mathit{_b}=1/56+1/146\,_{\textcolor{red}{ 2}}+1/584\,_{\textcolor{red}{ 8}}+1/1022\,_{\textcolor{red}{14}}}$ \\ \hline
$23$ & $47$ & $20$ & $15$ & $12$ &  $\mathbf{\textcolor{red}{ 3}}$& $60$ & $  \mathbf {2/73\mathit{_c}=1/60+1/219\,_{\textcolor{red}{ 3}}+1/292\,_{\textcolor{red}{ 4}}+1/365\,_{\textcolor{red}{5}}}\;\, ^{Eg}$ \\ \hline 
$8$ & $17$ & $12$ & $4$ & $1$ & $\mathbf{\textcolor{red}{ 3}}$& $48$ & $ \mathbf {2/79\mathit{_a}=1/48+1/316\,_{\textcolor{red}{ 4}}+1/948\,_{\textcolor{red}{ 12}}+1/3792\,_{\textcolor{red}{ 48}}}$ \\ \hline
$8$ & $17$ & $8$ & $6$ & $3$ & $\mathbf{\textcolor{red}{ 3}}$& $48$ & $ \mathbf {2/79\mathit{_b}=1/48+1/474\,_{\textcolor{red}{ 6}}+1/632\,_{\textcolor{red}{ 8}}+1/1264\,_{\textcolor{red}{ 16}}}$ \\ \hline
$16$ & $33$ & $28$ & $4$ & $1$ & $\mathbf{\textcolor{red}{ 3}}$& $56$ & $ \mathbf {2/79\mathit{_c}=1/56+1/158\,_{\textcolor{red}{ 2}}+1/1106\,_{\textcolor{red}{ 14}}+1/4424\,_{\textcolor{red}{ 56}}}$ \\ \hline 
$6$ & $13$ & $8$ & $4$ & $1$ & $\mathbf{\textcolor{red}{ 3}}$& $48$ & $ \mathbf {2/83\mathit{_a}=1/48+1/498\,_{\textcolor{red}{ 6}}+1/996\,_{\textcolor{red}{ 12}}+1/3984\,_{\textcolor{red}{ 48}}}$ \\ \hline
$8$ & $17$ & $10$ & $5$ & $2$ & $\mathbf{\textcolor{red}{ 3}}$& $50$ & $ \mathbf {2/83\mathit{_b}=1/50+1/415\,_{\textcolor{red}{ 5}}+1/830\,_{\textcolor{red}{ 10}}+1/2075\,_{\textcolor{red}{ 25}}}$ \\ \hline
$18$ & $37$ & $30$ & $5$ & $2$ & $\mathbf{\textcolor{red}{ 3}}$& $60$ & $ \mathbf {2/83\mathit{_c}=1/60+1/166\,_{\textcolor{red}{ 2}}+1/996\,_{\textcolor{red}{ 12}}+1/2490\,_{\textcolor{red}{30}}}$ \\ \hline \hline
$15$ & $31$ & $15$ & $10$ & $6$ & $\mathbf{\textcolor{red}{ 4}}$& $30$ & $ \mathbf {2/29=1/30+1/58\,_{\textcolor{red}{ 2}}+1/87\,_{\textcolor{red}{ 3}}+1/145\,_{\textcolor{red}{5}}}$ \\ \hline 
$14$ & $29$ & $15$ & $9$ & $5$ &  $\mathbf{\textcolor{red}{ 4}}$& $45$ & $ \mathbf {2/61=1/45+1/183\,_{\textcolor{red}{ 3}}+1/305\,_{\textcolor{red}{ 5}}+1/549\,_{\textcolor{red}{9}}}$ \\ \hline 
$17$ & $35$ & $27$ & $6$ & $2$ &  $\mathbf{\textcolor{red}{ 4}}$& $54$ & $ \mathbf {2/73=1/54+1/146\,_{\textcolor{red}{ 2}}+1/657\,_{\textcolor{red}{ 9}}+1/1971\,_{\textcolor{red}{ 27}}}$ \\ \hline
$15$ & $31$ & $15$ & $10$ & $6$ &  $\mathbf{\textcolor{red}{ 4}}$& $60$ & $ \mathbf {2/89=1/60+1/356\,_{\textcolor{red}{ 4}}+1/534\,_{\textcolor{red}{ 6}}+1/890\,_{\textcolor{red}{10}}}\;\, ^{Eg}$ \\ \hline \hline 
$9$ & $19$ & $12$ & $6  $ & $1 $ & $\mathbf{\textcolor{red}{ 5}}$& $24$ & $\mathbf {2/29=1/24+1/58\,_{\textcolor{red}{ 2}}+1/116\,_{\textcolor{red}{ 4}}+1/696\,_{\textcolor{red}{24}}}$ \\ \hline
$8$ & $17$ & $10$ & $6$ & $1$ & $\mathbf{\textcolor{red}{ 5}}$& $30$ & $ \mathbf {2/43=1/30+1/129\,_{\textcolor{red}{ 3}}+1/215\,_{\textcolor{red}{ 5}}+1/1290\,_{\textcolor{red}{30}}}$ \\ \hline
$11$ & $23$ & $14$ & $7$ & $2$ &  $\mathbf{\textcolor{red}{ 5}}$& $42$ & $ \mathbf {2/61\mathit{_a}=1/42+1/183\,_{\textcolor{red}{ 3}}+1/366\,_{\textcolor{red}{ 6}}+1/1281\,_{\textcolor{red}{21}}}$ \\ \hline
$17$ & $35$ & $24$ & $8$ & $3$ &  $\mathbf{\textcolor{red}{ 5}}$& $48$ & $ \mathbf {2/61\mathit{_b}=1/48+1/122\,_{\textcolor{red}{ 2}}+1/366\,_{\textcolor{red}{ 6}}+1/976\,_{\textcolor{red}{16}}}$ \\ \hline
$11$ & $23$ & $16$ & $6$ & $1$ &  $\mathbf{\textcolor{red}{ 5}}$& $48$ & $ \mathbf {2/73\mathit{_a}=1/48+1/219\,_{\textcolor{red}{ 3}}+1/584\,_{\textcolor{red}{ 8}}+1/3504\,_{\textcolor{red}{ 48}}}$ \\ \hline
$11$ & $23$ & $12$ & $8$ & $3$ &  $\mathbf{\textcolor{red}{ 5}}$& $48$ & $ \mathbf {2/73\mathit{_b}=1/48+1/292\,_{\textcolor{red}{ 4}}+1/438\,_{\textcolor{red}{ 6}}+1/1168\,_{\textcolor{red}{ 16}}}$ \\ \hline
$12$ & $25$ & $18$ & $6$ & $1$ & $\mathbf{\textcolor{red}{ 5}}$& $54$ & $ \mathbf {2/83\mathit{_a}=1/54+1/249\,_{\textcolor{red}{ 3}}+1/747\,_{\textcolor{red}{ 9}}+1/4482\,_{\textcolor{red}{ 54}}}$ \\ \hline
$18$ & $37$ & $30$ & $6$ & $1$ & $\mathbf{\textcolor{red}{ 5}}$& $60$ & $ \mathbf {2/83\mathit{_b}=1/60+1/166\,_{\textcolor{red}{ 2}}+1/830\,_{\textcolor{red}{ 10}}+1/4980\,_{\textcolor{red}{60}}}$ \\ \hline
$11$ & $23$ & $14$ & $7$ & $2$ &  $\mathbf{\textcolor{red}{ 5}}$& $56$ & $ \mathbf {2/89=1/56+1/356\,_{\textcolor{red}{ 4}}+1/712\,_{\textcolor{red}{ 8}}+1/2492\,_{\textcolor{red}{ 28}}}$ \\ \hline \hline
$14$ & $29$ & $18$ & $9$ & $2$ & $\mathbf{\textcolor{red}{ 7}}$& $36$ & $ \mathbf {2/43=1/36+1/86\,_{\textcolor{red}{ 2}}+1/172\,_{\textcolor{red}{ 4}}+1/774\,_{\textcolor{red}{18}}}$ \\ \hline
$9$ & $19$ & $10$ & $8$ & $1$ &  $\mathbf{\textcolor{red}{ 7}}$& $40$ & $\mathbf {2/61=1/40+1/244\,_{\textcolor{red}{ 4}}+1/305\,_{\textcolor{red}{ 5}}+1/2440\,_{\textcolor{red}{ 40}}}$ \\ \hline
$23$ & $47$ & $30$ & $12$ & $5$ &  $\mathbf{\textcolor{red}{ 7}}$& $60$ & $ \mathbf {2/73=1/60+1/146\,_{\textcolor{red}{ 2}}+1/365\,_{\textcolor{red}{ 5}}+1/876\,_{\textcolor{red}{12}}}$ \\ \hline 
$14$ & $29$ & $18$ & $9$ & $2$ & $\mathbf{\textcolor{red}{ 7}}$& $54$ & $ \mathbf {2/79=1/54+1/237\,_{\textcolor{red}{ 3}}+1/474\,_{\textcolor{red}{ 6}}+1/2133\,_{\textcolor{red}{ 27}}}$ \\ \hline
$18$ & $37$ & $20$ & $12$ & $5$ & $\mathbf{\textcolor{red}{ 7}}$& $60$ & $ \mathbf {2/83=1/60+1/249\,_{\textcolor{red}{ 3}}+1/415\,_{\textcolor{red}{ 5}}+1/996\,_{\textcolor{red}{12}}}$ \\ \hline 
$11$ & $23$ & $14$ & $8$ & $1$ &  $\mathbf{\textcolor{red}{ 7}}$& $56$ & $ \mathbf {2/89=1/56+1/356\,_{\textcolor{red}{ 4}}+1/623\,_{\textcolor{red}{ 7}}+1/4984\,_{\textcolor{red}{ 56}}}$ \\ \hline \hline
$20$ & $41$ & $21$ & $14$ & $6$ & $\mathbf{\textcolor{red}{ 8}}$& $42$ & $ \mathbf {2/43=1/42+1/86\,_{\textcolor{red}{ 2}}+1/129\,_{\textcolor{red}{ 3}}+1/301\,_{\textcolor{red}{7}}}\;\, ^{Eg}$ \\ \hline 
$15$ & $31$ & $15$ & $12$ & $4$ &  $\mathbf{\textcolor{red}{ 8}}$& $60$ & $ \mathbf {2/89=1/60+1/356\,_{\textcolor{red}{ 4}}+1/445\,_{\textcolor{red}{ 5}}+1/1335\,_{\textcolor{red}{ 15}}}$ \\ \hline \hline 
$21$ & $43$ & $26$ & $13$ & $4$ &  $\mathbf{\textcolor{red}{ 9}}$& $52$ & $ \mathbf {2/61=1/52+1/122\,_{\textcolor{red}{ 2}}+1/244\,_{\textcolor{red}{ 4}}+1/793\,_{\textcolor{red}{13}}}$ \\ \hline
$20$ & $41$ & $30$ & $10$ & $1$ & $\mathbf{\textcolor{red}{ 9}}$& $60$ & $ \mathbf {2/79\mathit{_a}=1/60+1/158\,_{\textcolor{red}{ 2}}+1/474\,_{\textcolor{red}{ 6}}+1/4740\,_{\textcolor{red}{ 60}}}$ \\ \hline
$20$ & $41$ & $20$ & $15$ & $6$ & $\mathbf{\textcolor{red}{ 9}}$& $60$ & $ \mathbf {2/79\mathit{_b}=1/60+1/237\,_{\textcolor{red}{ 3}}+1/316\,_{\textcolor{red}{ 4}}+1/790\,_{\textcolor{red}{10}}}\;\, ^{Eg}$ \\ \hline
$15$ & $31$ & $20$ & $10$ & $1$ &  $\mathbf{\textcolor{red}{ 9}}$& $60$ & $ \mathbf {2/89=1/60+1/267\,_{\textcolor{red}{ 3}}+1/534\,_{\textcolor{red}{ 6}}+1/5340\,_{\textcolor{red}{ 60}}}$ \\ \hline \hline
$25$ & $51$ & $35$ & $14$ & $2$ &  $\mathbf{\textcolor{red}{ 12}}$& $70$ & $ \mathbf {2/89=1/70+1/178\,_{\textcolor{red}{ 2}}+1/445\,_{\textcolor{red}{ 5}}+1/3115\,_{\textcolor{red}{35}}}$ \\ \hline \hline 
$23$ & $47$ & $30$ & $15$ & $2$ &  $\mathbf{\textcolor{red}{13}}$& $60$ & $ \mathbf {2/73=1/60+1/146\,_{\textcolor{red}{ 2}}+1/292\,_{\textcolor{red}{ 4}}+1/2190\,_{\textcolor{red}{30}}}$ \\ \hline
$18$ & $37$ & $20$ & $15$ & $2$ & $\mathbf{\textcolor{red}{ 13}}$& $60$ & $ \mathbf {2/83=1/60+1/249\,_{\textcolor{red}{ 3}}+1/332\,_{\textcolor{red}{ 4}}+1/2490\,_{\textcolor{red}{30}}}$ \\ \hline \hline
$24$ & $49$ & $32$ & $16$ & $1$ & $\mathbf{\textcolor{red}{ 15}}$& $64$ & $ \mathbf {2/79=1/64+1/158\,_{\textcolor{red}{ 2}}+1/316\,_{\textcolor{red}{ 4}}+1/5056\,_{\textcolor{red}{ 64}}}$ \\ \hline
$26$ & $53$ & $34$ & $17$ & $2$ & $\mathbf{\textcolor{red}{ 15}}$& $68$ & $ \mathbf {2/83=1/68+1/166\,_{\textcolor{red}{ 2}}+1/332\,_{\textcolor{red}{ 4}}+1/2822\,_{\textcolor{red}{34}}}$ \\ \hline \hline 
$23$ & $47$ & $27$ & $18$ & $2$ &  $\mathbf{\textcolor{red}{ 16}}$& $54$ & $ \mathbf {2/61=1/54+1/122\,_{\textcolor{red}{ 2}}+1/183\,_{\textcolor{red}{ 3}}+1/1647\,_{\textcolor{red}{27}}}$ \\ \hline \hline
$27$ & $55$ & $36$ & $18$ & $1$ &  $\mathbf{\textcolor{red}{ 17}}$& $72$ & $ \mathbf {2/89=1/72+1/178\,_{\textcolor{red}{ 2}}+1/356\,_{\textcolor{red}{ 4}}+1/6408\,_{\textcolor{red}{72}}}$ \\ \hline \hline
$30$ & $61$ & $36$ & $24$ & $1$ & $\mathbf{\textcolor{red}{ 23}}$& $72$ & $ \mathbf {2/83=1/72+1/166\,_{\textcolor{red}{ 2}}+1/249\,_{\textcolor{red}{ 3}}+1/5976\,_{\textcolor{red}{72}}}$ \\ \hline \hline
$33$ & $67$ & $39$ & $26$ & $2$ &  $\mathbf{\textcolor{red}{ 24}}$& $78$ & $ \mathbf {2/89=1/78+1/178\,_{\textcolor{red}{ 2}}+1/267\,_{\textcolor{red}{ 3}}+1/3471\,_{\textcolor{red}{39}}}$ \\ \hline 
\end{tabular}
\end{center}
\label{Complete4Terms}
\end{table}
\normalsize
\clearpage

Table \ref{Complete4Terms} shown above is only as an indication for us  and, certainly, was not calculated in its entirety. 
{\sf 2/23} has been reported only for memory because it was solved at the end of Sect. \ref{ThreeTerms}.\\
With their experience related to 3-terms series, cut-off beyond $10$ has been applied by the scribes. Indeed all cases (here \sf 7\rm) may support this cut-off without any exception.  Table \ref{Complete4Terms} becomes: \\
\begin{table}[h]
\caption{\sf [4-terms] options}
\small
\begin{center}
\begin{tabular}{|l|c|l||l||l||c|l|l|} \hline
\multicolumn{8}{|c|}{\sf Trials [4-terms] ordered with $\Delta_{d}^{'}\nearrow$ showing where are the Egyptian options}\\ \hline
$n$ & $2n+1$ & $d_2$ & $d_3$ & $d_4$ &$\textcolor{red}{\Delta_{d}^{'}}$  &$D_1^n$ & Possible [4-terms] decompositions $\mathbf {\textcolor{red}{ m_4\leq 10}}$ \\   \hline \hline
$9$ & $19$ & $12$ & $4  $ & $3 $ & $\mathbf{\textcolor{red}{ 1}}$& $24$ & $\mathbf {2/29=1/24+1/58\,_{\textcolor{red}{ 2}}+1/174\,_{\textcolor{red}{ 6}}+1/232\,_{\textcolor{red}{8}}}\;\, ^{Eg}$ \\   \hline \hline
$9$ & $19$ & $10$ & $5$ & $4$ &  $\mathbf{\textcolor{red}{ 1}}$& $40$ & $\mathbf {2/61=1/40+1/244\,_{\textcolor{red}{ 4}}+1/488\,_{\textcolor{red}{ 8}}+1/610\,_{\textcolor{red}{ 10}}}\;\, ^{Eg}$ \\   \hline \hline
$14$ & $29$ & $14$ & $8$ & $7$ & $\mathbf{\textcolor{red}{ 1}}$& $56$ & $ \mathbf {2/83=1/56+1/332\,_{\textcolor{red}{ 4}}+1/581\,_{\textcolor{red}{ 7}}+1/664\,_{\textcolor{red}{ 8}}}$ \\   \hline \hline \hline
$5$ & $11$ & $5$ & $4$ & $2$ & $\mathbf{\textcolor{red}{ 2}}$& $20$ & $ \mathbf {2/29=1/20+1/116\,_{\textcolor{red}{ 4}}+1/145\,_{\textcolor{red}{ 5}}+1/290\,_{\textcolor{red}{ 10}}}$ \\   \hline \hline
$18$ & $37$ & $15$ & $12$ & $10$ & $\mathbf{\textcolor{red}{ 2}}$& $60$ & $ \mathbf {2/83=1/60+1/332\,_{\textcolor{red}{ 4}}+1/415\,_{\textcolor{red}{ 5}}+1/498\,_{\textcolor{red}{6}}}\;\, ^{Eg{\textcolor{red}{{\star }}}}$ \\   \hline \hline
$18$ & $37$ & $21$ & $9$ & $7$ &  $\mathbf{\textcolor{red}{ 2}}$& $63$ & $ \mathbf {2/89=1/63+1/267\,_{\textcolor{red}{ 3}}+1/623\,_{\textcolor{red}{ 7}}+1/801\,_{\textcolor{red}{9}}}$ \\   \hline \hline \hline
$8$ & $17$ & $10$ & $5$ & $2$ & $\mathbf{\textcolor{red}{ 3}}$& $20$ & $\mathbf {\cancel{2/23}=1/20+1/46\,_{\textcolor{red}{ 2}}+1/92\,_{\textcolor{red}{ 4}}+1/230\,_{\textcolor{red}{10}}}$ \\ \hline \hline
$23$ & $47$ & $20$ & $15$ & $12$ &  $\mathbf{\textcolor{red}{ 3}}$& $60$ & $  \mathbf {2/73=1/60+1/219\,_{\textcolor{red}{ 3}}+1/292\,_{\textcolor{red}{ 4}}+1/365\,_{\textcolor{red}{5}}}\;\, ^{Eg{\textcolor{red}{{\star }}}}$ \\   \hline \hline \hline
$15$ & $31$ & $15$ & $10$ & $6$ & $\mathbf{\textcolor{red}{ 4}}$& $30$ & $ \mathbf {2/29=1/30+1/58\,_{\textcolor{red}{ 2}}+1/87\,_{\textcolor{red}{ 3}}+1/145\,_{\textcolor{red}{5}}}$ \\   \hline \hline
$14$ & $29$ & $15$ & $9$ & $5$ &  $\mathbf{\textcolor{red}{ 4}}$& $45$ & $ \mathbf {2/61=1/45+1/183\,_{\textcolor{red}{ 3}}+1/305\,_{\textcolor{red}{ 5}}+1/549\,_{\textcolor{red}{9}}}$ \\   \hline \hline
$15$ & $31$ & $15$ & $10$ & $6$ &  $\mathbf{\textcolor{red}{ 4}}$& $60$ & $ \mathbf {2/89=1/60+1/356\,_{\textcolor{red}{ 4}}+1/534\,_{\textcolor{red}{ 6}}+1/890\,_{\textcolor{red}{10}}}\;\, ^{Eg}$ \\   \hline \hline \hline
$20$ & $41$ & $21$ & $14$ & $6$ & $\mathbf{\textcolor{red}{ 8}}$& $42$ & $ \mathbf {2/43=1/42+1/86\,_{\textcolor{red}{ 2}}+1/129\,_{\textcolor{red}{ 3}}+1/301\,_{\textcolor{red}{7}}}\;\, ^{Eg}$ \\   \hline \hline
$20$ & $41$ & $20$ & $15$ & $6$ & $\mathbf{\textcolor{red}{ 9}}$& $60$ & $ \mathbf {2/79=1/60+1/237\,_{\textcolor{red}{ 3}}+1/316\,_{\textcolor{red}{ 4}}+1/790\,_{\textcolor{red}{10}}}\;\, ^{Eg}$ \\   \hline 
\end{tabular}
\end{center}
\label{4TERMSOPT}
\end{table}
We follow the same way as for the [3-terms] series with slightly different subsets. That yields:\\
\begin{table}[htbp]
\caption{\sf A single or two different $\Delta_{d}^{'}$  [4-terms]\rm }
\begin{center}
\scriptsize
\begin{tabular}{|l|c|l||l||l||c|l|l|} \hline
\multicolumn{7}{|c|}{\sf D with a single $\Delta_{d}^{'} $$\quad$(options: no)} & \multicolumn{1}{l|}{\sf Scribes's decision: obvious }\\ \hline
$n$ & $2n+1$ & $d_2$ & $d_3$ & $d_4$ &$\textcolor{red}{\Delta_{d}^{'}}$ &$D_1^n$ &$\qquad \quad$[4-terms] decompositions  \\   \hline \hline
$23$ & $47$ & $20$ & $15$ & $12$ &  $\mathbf{\textcolor{red}{ 3}}$& $60$ & $  \mathbf {2/73=1/60+1/219\,_{\textcolor{red}{ 3}}+1/292\,_{\textcolor{red}{ 4}}+1/365\,_{\textcolor{red}{5}}}\;\, ^{Eg{\textcolor{red}{{\star }}}}$ \\   \hline \hline \hline
$20$ & $41$ & $21$ & $14$ & $6$ & $\mathbf{\textcolor{red}{ 8}}$& $42$ & $ \mathbf {2/43=1/42+1/86\,_{\textcolor{red}{ 2}}+1/129\,_{\textcolor{red}{ 3}}+1/301\,_{\textcolor{red}{7}}}\;\, ^{Eg}$ \\   \hline \hline \hline
$20$ & $41$ & $20$ & $15$ & $6$ & $\mathbf{\textcolor{red}{ 9}}$& $60$ & $ \mathbf {2/79=1/60+1/237\,_{\textcolor{red}{ 3}}+1/316\,_{\textcolor{red}{ 4}}+1/790\,_{\textcolor{red}{10}}}\;\, ^{Eg}$ \\   \hline 
\end{tabular}\\ \vspace{0.5em}
\begin{tabular}{|l|c|l||l||l||c|l|l|} \hline
\multicolumn{7}{|c|}{\sf D with two different  $\Delta_{d}^{'}\quad$(options: yes)} & \multicolumn{1}{l|}{\sf Scribes's decision: smallest  $\Delta_{d}^{'} $}\\ \hline
$n$ & $2n+1$ & $d_2$ & $d_3$ & $d_4$ &$\textcolor{red}{\Delta_{d}^{'}}$ &$D_1^n$ & $\qquad \quad$[4-terms] decompositions \\   \hline \hline
$9$ & $19$ & $12$ & $4  $ & $3 $ & $\mathbf{\textcolor{red}{ 1}}$& $24$ & $\mathbf {2/29=1/24+1/58\,_{\textcolor{red}{ 2}}+1/174\,_{\textcolor{red}{ 6}}+1/232\,_{\textcolor{red}{8}}}\;\, ^{Eg}$ \\   \hline \hline
$5$ & $11$ & $5$ & $4$ & $2$ & $\mathbf{\textcolor{red}{ 2}}$& $20$ & $ \mathbf {2/29=1/20+1/116\,_{\textcolor{red}{ 4}}+1/145\,_{\textcolor{red}{ 5}}+1/290\,_{\textcolor{red}{ 10}}}$ \\   \hline \hline
$15$ & $31$ & $15$ & $10$ & $6$ & $\mathbf{\textcolor{red}{ 4}}$& $30$ & $ \mathbf {2/29=1/30+1/58\,_{\textcolor{red}{ 2}}+1/87\,_{\textcolor{red}{ 3}}+1/145\,_{\textcolor{red}{5}}}$ \\   \hline \hline \hline
$9$ & $19$ & $10$ & $5$ & $4$ &  $\mathbf{\textcolor{red}{ 1}}$& $40$ & $\mathbf {2/61=1/40+1/244\,_{\textcolor{red}{ 4}}+1/488\,_{\textcolor{red}{ 8}}+1/610\,_{\textcolor{red}{ 10}}}\;\, ^{Eg}$ \\   \hline \hline
$14$ & $29$ & $15$ & $9$ & $5$ &  $\mathbf{\textcolor{red}{ 4}}$& $45$ & $ \mathbf {2/61=1/45+1/183\,_{\textcolor{red}{ 3}}+1/305\,_{\textcolor{red}{ 5}}+1/549\,_{\textcolor{red}{9}}}$ \\  \hline
\end{tabular}\\ \vspace{0.5em}
\begin{tabular}{|l|c|l||l||l||c|l|l|} \hline
\multicolumn{7}{|c|}{} & \multicolumn{1}{l|}{\sf Scribes's decision: consecutive multipliers}\\ \hline
$n$ & $2n+1$ & $d_2$ & $d_3$ & $d_4$ &$\textcolor{red}{\Delta_{d}^{'}}$&$D_1^n$ & $\qquad \quad$[4-terms] decompositions \\   \hline \hline
$14$ & $29$ & $14$ & $8$ & $7$ & $\mathbf{\textcolor{red}{ 1}}$& $56$ & $ \mathbf {2/83=1/56+1/332\,_{\textcolor{red}{ 4}}+1/581\,_{\textcolor{red}{ 7}}+1/664\,_{\textcolor{red}{ 8}}}$ \\   \hline \hline
$18$ & $37$ & $15$ & $12$ & $10$ & $\mathbf{\textcolor{red}{ 2}}$& $60$ & $ \mathbf {2/83=1/60+1/332\,_{\textcolor{red}{ 4}}+1/415\,_{\textcolor{red}{ 5}}+1/498\,_{\textcolor{red}{6}}}\;\, ^{Eg{\textcolor{red}{{\star }}}}$ \\   \hline 
\end{tabular}\\ \vspace{0.5em}
\begin{tabular}{|l|c|l||l||l||c|l|l|} \hline
\multicolumn{7}{|c|}{} & \multicolumn{1}{l|}{\sf Scribes's decision: no odd denominator $D_1$}\\ \hline
$n$ & $2n+1$ & $d_2$ & $d_3$ & $d_4$ &$\textcolor{red}{\Delta_{d}^{'}}$&$D_1^n$ & $\qquad \quad$[4-terms] decompositions \\   \hline \hline
$18$ & $37$ & $21$ & $9$ & $7$ &  $\mathbf{\textcolor{red}{ 2}}$& $63$ & $ \mathbf {2/89=1/63+1/267\,_{\textcolor{red}{ 3}}+1/623\,_{\textcolor{red}{ 7}}+1/801\,_{\textcolor{red}{9}}}$ \\   \hline \hline 
$15$ & $31$ & $15$ & $10$ & $6$ &  $\mathbf{\textcolor{red}{ 4}}$& $60$ & $ \mathbf {2/89=1/60+1/356\,_{\textcolor{red}{ 4}}+1/534\,_{\textcolor{red}{ 6}}+1/890\,_{\textcolor{red}{10}}}\;\, ^{Eg}$ \\   \hline
\end{tabular}	
\end{center}
\normalsize
\label{1Delta4}
\end{table}
We recall that any odd denominator $D_1$ could lead to a solution for  [3-terms] decompositions as checked in tables  \ref {3TERMSOPT} or  \ref  {Frac3become4}. Its occurrence arises only 2 times in table \ref{1Delta4} [4-terms]. The first, for $2/61$, 
was dropped out because a ${\Delta_{d}^{'}}=4$ too high. The second one regards {\sf 2/89} (first row). Then, for a unifying sake and avoiding singularity, chief scribe decided to discard $D_1=63$ in this case.\\ 
\hspace*{1.5em}Remark that we are  very far from assumptions of Gillings {\bf  \cite{Gillings}} 
about Egyptian preferences for even numbers instead of odd, regarding the denominators in general.
Thus the {\it `no odd precept' } was a low priority. At  low ratio also (2 times only), this will be applied to the composite numbers $D$ {\bf  \cite{Brehamet}}.

										\section {Conclusion}
 \hspace{1.5em}As we saw, the most recent  analysis (2008) has been performed  on the `$2/n$' table
by Abdulaziz {\bf   \cite{Abdulaziz}} (see his group $G_2$). It can be appreciated as a kind of mathematical anastylosis using materials issued from the RMP and other documentation.
Ancient calculation procedure, using mainly fractions, is faithfully respected, but leads to arithmetical depth analyses of each divisor of $D_1$. \\

\hspace*{1.5em}Our global approach avoided the difficulties of sophisticated arithmetical studies.
This provides the advantage of forgetting quickly some widespread `modern' ideas about the topic. \\
$\bullet$ No, the last denominator is not bounded  by a fixed value of $1000$.
It only depends on the `circumstances'  related to the value of $D$. For 3 or 4 terms, a limitation like $D_h \leq 10D$ is quite suitable, except only for {\sf 2/53} where $10$ is replaced by $15$. An observation well stressed in Ref. {\bf   \cite{Abdulaziz}}.\\
$\bullet$ No requirement is found about the  denominator $D_1$ as having to be the greatest if alternatives.\\
$\bullet$ Once for all, a systematic predilection for even denominators  does not need to be considered. Only once, we were forced to discard $D_1 = 63$ (odd) for deciding on {\sf 2/89}.\\
$\bullet$ Of course, there is no theoretical formula that can give immediately the first denominator as a function of $D$. It must necessarily go through trials and few selection criteria. The simpler the better, like the $\Delta$-classification presented in  this paper.
Maybe is it this classification that induces  the opportunity  of a comprehensive approach ? Strictly speaking, there are no algorithms in the method, just tables and pertinent observation. This is how $2/23$,  $2/29 $ or $2/53$ have found a logical explanation, more thorough than the arguments commonly supplied for these `singularities'. \\
\hspace*{1.5em}Find a simple logic according to which there is no singular case was the goal of the present paper.\\
Perhaps, chronologically, the study of prime numbers has been elaborated \sf before \rm that of composite numbers. It is nothing more than an hypothesis consistent with the spirit of our study. 
Yes ancient scribes certainly have been able to calculate and analyze all the preliminary cases.
Ultimately, our unconventional method allows to reconstruct the table fairly easily with weak mathematical assumptions, except maybe the new idea to consider as beneficial to have consecutive multipliers.

\renewcommand{\theequation}{A.\arabic{equation}}
\setcounter{equation}{0}
\section*{\sf Appendix A: why a boundary with a Top-flag?}
In this appendix, we continue to  consider prime denominators $D$. For [2-terms] decompositions this concept of a Top-flag has no meaning since the last denominator is unique. \\
Obviously, doubtless far from Egyptian concepts, there are another equations {\sf  more general} than\\ 
Eqs. (\ref{eq:FEgypt3}) or (\ref{eq:FEgypt4}), namely
\begin{equation}
\frac{2}{D}= \frac{1}{D_1}+ \frac{1}{m_2D}+ \frac{1}{m_3D}.
\end{equation}
\begin{equation}
\frac{2}{D}= \frac{1}{D_1}+ \frac{1}{m_2D}+ \frac{1}{m_3D}+ \frac{1}{m_4D}.
\end{equation}
We can imagine these as issued from another kind of unity decomposition like
\begin{equation}
\mathbf{1}= \frac{D}{2D_1}+ \frac{1}{2m_2}+  \frac{1}{2m_3}.
\end{equation}
\begin{equation}
\mathbf{1}= \frac{D}{2D_1}+ \frac{1}{2m_2}+  \frac{1}{2m_3}+  \frac{1}{2m_4}.
\end{equation}
$D/2D_1$ remains in the lead of equality and $\mathbf{1}$ is a sum of terms, each with a even denominator . \\
These (modern) equations have additional solutions of no use for the scribes .\\

\it A priori \rm the solutions are infinite, then for avoiding such a tedious research (today and in the past time), it is necessary to limit
the highest denominator $D_h=m_h D$. How to do that ?
Simply by defining a kind of `Top-flag' $\mbox{\boldmath $\top$ }\!\!_ f^{\;\;[h]}$ such as
\begin{equation}
D_h \leq D \mbox{\boldmath $\top$ }\!\!_ f^{\;\;[h]}  .
\end{equation}
Indeed, as soon as one decides to study a three-terms decomposition or more,  it should be realized that an upper boundary for the last denominator has to be fixed. 
If not, the number of solutions becomes infinite [countable]. Recall that $m_2 < m_3 < m_4$ and $D_2 < D_3 < D_4$.
Unfortunately (or not)  the author of this paper has begun the calculations with a even more general problem, this of solving 
\begin{equation}
\frac{2}{D}= \sum _{i=1}^{h} \frac{1}{D_i},
\label{eq:EgyptGeneral}
\end{equation}
without any criteria of multiplicity involving multipliers like $m_i$ ($i>2$).\\
Certainly this was the reflex of Gillings {\bf{\cite{Gillings}}} or Bruckheimer and Salomon {\bf{\cite{BruckSalom}}}.
The problem is solvable and the solutions available by means of a small computer. After a necessary arithmetical  analysis,
it can be found that $(h-1)$ sets of solutions exist. One with $(h-1)$ multipliers $m_i$, another  with $(h-2)$ multipliers and so on.
No solution exists if one searches for $D_i$ ($i\geq 2$) not multiple of $D$.\\
Even a  low-level programming code like {\sf sb} can  be used instead of  {\sf Fortran}
to perform computations in a very acceptable speed. 
We  quickly realized the necessity of  stopping the calculations by using a limitation regarding the last highest denominator $D_h$. 
Whence the introduction of a Top-flag.\\
Actually the Egyptian $2/D$ table shows  a subset of more general solutions because the multipliers $m_i$ have a specific form 
involving $D_\mathbf1$ and some of its divisors $d_i$.
For example out of this subset, you can find an unexpected [4-terms] solution for 2/23 with $\mbox{\boldmath $\top$ }\!\!_ f^{\;\;[4]}= 10$, namely
\begin{tabular}{|c|}
\hline
$2/23=1/15+1/115\,_{\textcolor{red}{ 5}}+1/138\,_{\textcolor{red}{ 6}}+1/230\,_{\textcolor{red}{ 10}}$   \\ [0.01in] \hline
\end{tabular}.\\
\hspace*{1.5em}So, if we restrict ourself to retrieve Egyptian fractions given in the table, it naturally comes to mind to limit
the highest denominator by an upper boundary: a convenient Top-flag. \\
Excepted the Babylonian system example in base $60$,  a numeration in base $10$ is rather universal, because of our two hands with each {\sf 5} fingers.  
It is of common sense that the selection was generally $\mbox{\boldmath $\top$ }\!\!_ f^{\;\;[h]}= 10 \;(=2\times \sf 5)$, not excluding a favorable appreciation for 
$\mbox{\boldmath $\top$ }\!\!_ f^{\;\;[3]}= 15 \;(=3\times \sf 5)$ as for {\sf 2/53}.
\end{flushleft}



\begin{thebibliography}{99}


\bibitem {Peet} T. E. PEET: 
\em {The Rhind Mathematical Papyrus, British Museum 10057 and 10058}, \rm \\
London: The University Press of Liverpool limited and Hodder - Stoughton limited (1923).

\bibitem{Chace}A. B. CHACE; l. BULL; H. P. MANNING; and R. C. ARCHIBALD:
{\em The Rhind Mathematical Papyrus,} Mathematical Association of America, {\bf Vol.1 }(1927),
{\bf Vol. 2 }(1929), Oberlin, Ohio.

\bibitem{Robins} G. ROBINS and C. SHUTE:
{\em The Rhind Mathematical Papyrus: An Ancient Egyptian Text}\rm, London: British Museum Publications Limited,  (1987).
[A recent overview].

\bibitem{Gillings}R.J. GILLINGS: {\em Mathematics in the Time of Pharaohs}, MIT Press (1972), reprinted by Dover Publications (1982).

\bibitem{BruckSalom}M. BRUCKHEIMER  and Y. SALOMON: 
{\em Some comments on R.J Gillings's analysis of the 2/n table in the Rhind Papyrus}\rm , Historia Mathematica, {\bf Vol. 4}, pp. 445-452 (1977).


\bibitem{Imhausen} A. IMHAUSEN and J. RITTER:
\em{Mathematical fragments [see fragment UC32159]}. \rm (2004). \\ 
In: The UCL Lahun Papyri, {\bf Vol. 2} , pp. 71-96. 
Archeopress, Oxford, \\ Eds M. COLLIER, S. QUIRKE.

\bibitem{Abdulaziz}A. ABDULAZIZ: 
\em{ On the Egyptian method  of decomposing 2/n into unit fractions}\rm , Historia Mathematica, {\bf Vol. 35}, pp. 1-18 (2008).

\bibitem{GardnerMilo} M. GARDNER: {\em Egyptian fractions:} Unit Fractions, Hekats and Wages - an Update (2013), available on the site of academia.edu.
[Herein can be found an historic of various researches about the subject]. 


\bibitem{Brehamet}L. BREHAMET: {\em Remarks on the Egyptian 2/D table in favor of a global approach
(D composite number)}\rm,  arXiv [math.HO], to be submitted.


\bibitem {Fibonacci} L. FIBONACCI:  {\em Liber abaci} (1202).

\bibitem{Bruins}E.M. BRUINS: {\em The part in ancient Egyptian mathematics}, Centaurus, {\bf Vol. 19}, pp. 241-251 (1975).






 

				\end{thebibliography}
\end{document}